\documentclass[reqno,11pt]{amsart}

\usepackage{fbb} % osf for text, lining for math
\usepackage[varqu,varl]{inconsolata}% typewriter
\usepackage[libertine,bigdelims,vvarbb]{newtxmath}

\usepackage[cal=boondoxo]{mathalfa}% less slanted than STIX
\usepackage[T1]{fontenc}
\useosf

\usepackage{bbold}

\usepackage{stmaryrd}
\usepackage[mathscr]{eucal}

\usepackage{tikz-cd}

\usepackage{enumerate}
\usepackage[all]{xy}

\providecommand{\U}[1]{\protect\rule{.1in}{.1in}}
\newtheorem{theorem}{Theorem}
\theoremstyle{plain}

\newtheorem{claim}[theorem]{Claim}
\newtheorem{corollary}[theorem]{Corollary}
\newtheorem{lemma}[theorem]{Lemma}

\newtheorem{proposition}[theorem]{Proposition}

\theoremstyle{definition}
\newtheorem{example}[theorem]{Example}
\newtheorem{notation}[theorem]{Notation}
\newtheorem{remark}[theorem]{Remark}
\newtheorem{definition}[theorem]{Definition}

\numberwithin{equation}{section}
\numberwithin{theorem}{section}

\newcommand{\perHA}{\mathfrak{S}\mathsf{Sym}}
\newcommand{\twopii}{(2 \pi \mbi)}

\newcommand{\reg}{\mathsf{Reg}}
\newcommand{\dcz}[3]{H^{#1}_{\scrD}(#3;\bbz(#2))} %  Deligne coh. (for complex varieties) w/coeffs. in Z
\newcommand{\dcx}[3]{\bbz(#1)_\scrD^{#2}(#3)} % Deligne complex of currents computing DC for qp vars.

\newcommand{\llb}{\llbracket}
\newcommand{\rrb}{\rrbracket}
\newcommand{\la}{\langle}
\newcommand{\ra}{\rangle}
\newcommand{\norm}[1]{|\hspace{-1pt}| #1 |\hspace{-1pt}| }

\newcommand{\Hom}[3]{\mathsf{Hom}_{#1}(#2,#3)}
\newcommand{\zhom}[3]{\bbz\mathsf{Hom}_{#1}(#2,#3)}

\newcommand{\asimplex}[1]{\Delta^{#1}}
\newcommand{\acube}[1]{\oblong^{#1}}

\newcommand{\uni}[1]{\mathbb{\Gamma}^{#1}}
\newcommand{\unip}[2]{\mathbb{\Gamma}_{\hspace{-0.1cm}\Pi}^{#1}(#2)}

\newcommand{\cur}[2]{{'\mathscr{D}}^{#1}(#2)}    % Currents of degree  {#1}  on {#2}. 
\newcommand{\curd}[2]{{'\mathscr{D}}_{#1}(#2)}    % Currents of dimension  {#1}  on {#2}. 

\newcommand{\norlog}[3]{{\mathscr{N}}^{#1}(#2)\la \log{#3} \ra} % Normal Currents of degree  #1 on #2 with log poles on #3
\newcommand{\lnor}[2]{\mathscr{N}_{\mathrm{loc}}^{#1}(#2)} % Locally normal currents of degree #1 on #2.
\newcommand{\nor}[2]{\mathscr{N}^{#1}(#2)} % Normal currents of degree #1 on #2.
\newcommand{\intcur}[2]{\mathscr{I}^{#1}(#2)} % Integral currents of degree #1 on #2.

\newcommand{\dif}[1]{\mathsf{dif}_{\!  #1}}
\newcommand{\symg}[1]{\mathfrak{S}_{#1}}

\newcommand{\shuf}[2]{\mathsf{Sh}(#1, #2)}

\newcommand{\sch}[1]{\mathsf{Sch}_{#1}}

\newcommand{\btr}[1]{\blacktriangle^{\hspace{-0.04cm} #1}}

\newcommand{\vtriang}{
\Delta
  \kern-0.56em\mathord{\raisebox{-0.5\depth}{\scalebox{0.64}{\( \Delta\)} } } 
  \kern-0.616em\mathord{\raisebox{-0.5\depth}{\scalebox{0.36}{\( \Delta\)} } } 
 }
\newcommand{\simplex}[1]{\vtriang^{\hspace{-1pt}#1}} %% normal
\newcommand{\om}[2]{\omega_{#1}^{#2}}

\newcommand{\ve}{\varepsilon}
\newcommand{\hboxt}{{\hat \boxtimes_t}}

\newcommand{\tet}[2]{\theta_{#1}^{#2}}

\newcommand{\sbcxeq}[3]{\mathcal{Z}_{\text{eq}}^{#1}(#3; #2)}% simplicial Bloch equidimensional HCG complex
\newcommand{\shcg}[3]{CH^{#1}(#3; #2)}  % simplicial HCG
\newcommand{\msfF}{\mathsf{F}}

\newcommand{\msfT}{\mathsf{T}}
\newcommand{\grF}[1]{\mathsf{gr}^{#1}_\msfF} 

 \newcommand{\bone}{\mathbb{1}}
\newcommand{\mba}{\mathbf{a}}
 \newcommand{\mbb}{\mathbf{b}}
\newcommand{\mbd}{\mathbf{d}}
 \newcommand{\mbe}{\mathbf{e}}
 \newcommand{\mbi}{{\rm{i}}}
\newcommand{\mbu}{\mathbf{u}}
 \newcommand{\mbv}{\mathbf{v}}
 \newcommand{\mbt}{\mathbf{t}}
 \newcommand{\mbx}{\mathbf{x}}
\newcommand{\mby}{\mathbf{y}}
 \newcommand{\mbz}{\mathbf{z}}
\newcommand{\bba}{\mathbb{A}}
\newcommand{\bbc}{\mathbb{C}}
\newcommand{\bbf}{\mathbb{F}}

\newcommand{\bbn}{\mathbb{N}}
\newcommand{\bbp}{\mathbb{P}}
\newcommand{\bbq}{\mathbb{Q}}
\newcommand{\bbr}{\mathbb{R}}
\newcommand{\bbw}{\mathbb{W}}
\newcommand{\bbz}{\mathbb{Z}}

\newcommand{\calc}{\mathcal{C}}
\newcommand{\cald}{\mathcal{D}}
\newcommand{\call}{\mathcal{L}}
\newcommand{\cals}{\mathcal{S}}
\newcommand{\scrA}{\mathscr{A}}
\newcommand{\scrC}{\mathscr{C}}
\newcommand{\scrD}{\mathscr{D}}
\newcommand{\scrF}{\mathscr{F}}
\newcommand{\scrI}{\mathscr{I}}
\newcommand{\scrM}{\mathscr{M}}
\newcommand{\scrN}{\mathscr{N}}

\newcommand{\scrQ}{\mathscr{Q}}
\newcommand{\msfP}{\mathsf{P}}

\newcommand{\balpha}{\boldsymbol{\alpha}}
\newcommand{\bbeta}{\boldsymbol{\beta}}

\newcommand{\bchi}{\boldsymbol{\chi}}

\allowdisplaybreaks

\title[Multiplicative properties of regulators]{Multiplicative properties of the current transform regulator}

\author[Paulo Lima-Filho]{Paulo Lima-Filho}
\address{Department of Mathematics, Texas A\&M University, College Station, TX 77843}
\email{plfilho@math.tamu.edu}
\thanks{The author thanks the support and hospitality of the Isaac Newton Institute for Mathematical Sciences, during the program "K-theory, algebraic cycles and motivic homotopy theory" (KAH2, 2022), where this work was concludedx.}
\dedicatory{}

\begin{document}

\begin{abstract}
%\lipsum[3]
This paper utilizes the properties of transforms of currents under equidimensional cycles, as introduced in \cite{MR4498559}, to establish the multiplicative nature of the resulting regulator map, in the derived category. The construction relies on a synthetic presentation of the fundamental triples of currents from \cite{MR4498559}, which exhibits group-like behavior under an extended Eilenberg-Zilber morphism. A key component of the analysis is a character of the permutation Hopf algebra  \( \perHA\) that takes values in the function field of \( \bba_{\bbq}^\infty\).
\end{abstract}

\maketitle
\tableofcontents
%

%%%%%%%%%%%%. 

%\input{Intro_final.tex}
% !TEX root =  Lima-Filho-HCG_Reg_Prod-501 copy.tex

\section*{Introduction}

The study of regulator maps \( \reg \colon \shcg{p}{n}{U}\to \dcz{2p-n}{p}{U}\)  from the higher Chow groups of quasiprojective varieties  to their Deligne-Beilinson cohomology has been the subject of extensive research, and approached from many distinct perspectives. Explicit formulas for such maps can be found, for example, in  \cite{MR2342640}, \cite{MR3722692} and \cite{MR4498559}. In this paper we employ  \emph{current transform} techniques,  introduced in \cite{MR4498559}, to study  the multiplicative properties of the corresponding regulator maps.

From the start we use Suslin's equidimensional complexes \( \sbcxeq{p}{*}{X}\) to define higher Chow groups \cite{Sus-HighCh}. 
The advantage of representing a higher Chow class by a cycle \( \Upsilon\) in \( U\times \Delta^n\) which is equidimensional and dominant cycle over the algebraic simplex \( \Delta^n\),  is that \( \Upsilon \) induces a transform operation \(\Upsilon^\vee\) from suitable classes of currents in the compactification \( \bbp^n\) of \( \Delta^n\) to currents in the variety \( U\); see \eqref{eq:curr_transf}.

A crucial ingredient in our formalism is a triple of currents 
\( \bbf_n = (\simplex{n}, \Theta_n, W_n) \) in complex projective space \( \bbp^n\); see \S \ref{subsec:FTC}.   When 
a class \( \alpha \in \shcg{p}{n}{U}\) is represented by an equidimensional cycle \( \Upsilon \) in  \( U \times \Delta^n\), 
the transform of \( \bbf_n \) under \( \Upsilon\) 
gives a triple of currents 
\( \Upsilon^\vee(\bbf_n) = ( \Upsilon^\vee \simplex{n}, \Upsilon^\vee \Theta_n, \Upsilon^\vee W_n) \) on \( U \) that  represents \( \reg (\alpha)  \in \dcz{2p-n}{p}{U}\).

In reality, this construction gives a map of complexes
\[
\reg_U \ \colon \ \sbcxeq{p}{*}{U} \to \bbz(p)^{2p-*}_\scrD(U),
\]
where \(\bbz(p)^{*}_\scrD(U)\)  is a complex of currents whose cohomology is 
\( \dcz{*}{p}{U}\), 
%Despite utilizing a good compactification
%\( 
%\begin{tikzcd}[column sep=small]
%U \ar[r, hook] & X & D = X-U \ar[l, hook']
%\end{tikzcd}
%\)
%of \( U\), 
 which becomes a natural transformation of functors into the derived category of bounded complexes of abelian groups.

Remarkably, the triples $\bbf_n$ possess highly structured properties that can be elegantly packaged and, in conjunction with the current transforms formalism, provide a  conceptually simple proof of the following result.
\medskip

\noindent{\bf Theorem \ref{thm:MAIN}.}{\it \ \ 
Let \( U \) and \( V\) be smooth quasi-projective varieties over \( \bbc\). Then, for any fixed \( t \in \bbc\), the following diagram of morphism of complexes commutes up to natural homotopy
\begin{equation}
\begin{tikzcd}
\sbcxeq{p}{*}{U} \otimes \sbcxeq{q}{*}{V} \ar[rr,"\bar \times"] \ar[d, "\reg_U\otimes\,  \reg_V"']  & &  \sbcxeq{p+q}{*}{U\times V} 
\ar[d, "\reg_{U\times V}"] \\
\dcx{p}{2p-*}{U}\otimes \dcx{q}{2q-*}{V} \ar[rr,"\boxtimes_t"'] & &  \dcx{p+q}{2(p+q)-*}{U\times V},
\end{tikzcd}
\end{equation}
where 
\( \bar \times \) and \( \boxtimes_t \) are the exterior products in the equidimensional higher Chow complexes and Deligne-Beilinson current complexes, respectively.
}
\medskip

The expected consequences follow from this theorem. 
\medskip

\noindent{\bf Corollary \ref{cor:ext}.}{\it \ \ 
If \( U \) and \( V\) are smooth quasi-projective varieties over \( \bbc\), then the regulator maps 
%\[ \reg_U \colon \shcg{p}{r}{U}\to \dcz{2p-r}{p}{U}\] and \( \reg_V \colon \shcg{q}{s}{V}\to \dcz{2q-s}{q}{V}\) 
commute with exterior product. In other words,
\[ \reg_{U\times V}(\alpha \times_{\text{\tiny CH}} \beta)  = \reg_U(\alpha) \times_\scrD \reg_V(\beta)\  \in\  \dcz{2(p+q)-(r+s)}{p+q}{U\times V},\] where 
\(  \times_{\text{\tiny CH}} \) and \( \times_\scrD\) are the exterior products in higher Chow groups and Deligne-Beilinson cohomology, respectively.
}
\medskip

Using  the natural isomorphism   \( \shcg{p}{n}{U}\cong H_\scrM^{2p-n}(U;\bbz(p)) \) between higher Chow groups and Voevodsky's motivic cohomology  \cite[Lect. 19]{MR2242284}, one obtains the next result. 
\medskip

\noindent{\bf Corollary \ref{cor:cup}.}{\it \ \ 
The regulator map \( \reg_U \colon  H_\scrM^{*}(U;\bbz(\bullet)) \to  \dcz{*}{\bullet}{U}\) is a homomorphism of bigraded rings, and is a natural transformation of functors in the category of smooth complex quasiprojective varieties.
}
\medskip

To describe the main arguments,  we need to introduce a graded chain complex  \(\,  (\unip{*}{\bullet}, D) \), as follows.  
For each \( k\)-tuple of non-negative integers \( \mbb = (b_1, \ldots, b_k) \), denote  \( |\mbb| = \sum_{j=1}^k b_j \)  and set  \( \bbp(\mbb) = \bbp^{b_1} \times \cdots \times \bbp^{b_k} \). 
Then define
\begin{equation}
\unip{r}{k} \ := \ \prod_{\mbb \in \bbz_{\geq 0}^k} \cur{|\mbb|  + r}{  \bbp(\mbb)},
\end{equation}
where   \( \cur{d}{X} \) denotes the group of deRham currents  of degree \(d \) on a manifold \( X\). The differential \( D \colon \unip{r}{k} \to \unip{r+1}{k} \) is induced by the exterior derivative of currents along with the push-forward of currents under inclusions of coordinate subspaces. In reality, \( (\unip{*}{k}, D)\) is the total product complex of a \((k+1)\)-multicomplex; as explained in \S\ref{subsec:univ}. 

While we work with currents from the outset, implicit in the background are  \emph{multi semi-cosimplicial varieties} (no degeneracies) constructed from  the projective compactifications \( \bbp(\mbb)\) of products of algebraic simplices \( \Delta^{b_1}\times \cdots \times \Delta^{b_k}\).  In this context, our constructions work because  currents form differential co-sheaves in the category of smooth manifolds and proper maps.

The Hodge filtration on currents induces a filtration 
\[ \cdots \subset F^{-1} \unip{*}{k} \subset  F^{0} \unip{*}{k} \subset  F^{1} \unip{*}{k} \subset \cdots \subset  F^{p} \unip{*}{k} \subset \cdots \] 
 on the complex \( \unip{*}{k} \)
extending in both directions, which restricts to corresponding filtrations on  several subcomplexes of interest. Particularly,  the subcomplexes \( \scrI_\Pi^*(k) \subset \scrN^*_\Pi(k) \subset \unip{*}{k}\), where \( \scrI^* \) and \( \scrN^* \) denote  integral and normal currents, respectively, and  the associated subcomplexes \( \scrI_\Pi^*(k|\infty) \subset \scrN^*_\Pi(k|\infty) \subset \unip{*}{k|\infty} \subset \unip{*}{k}\)  of currents that \emph{vanish suitably at infinity}; see \S\ref{subsec:subcxs}  and \cite[Defn. 3]{MR4498559}.

The main arguments can be fully appreciated once we simplify  the notation by writing   elements \(\alpha \in  \unip{r}{k} \) as   power series \( \alpha = \sum_{\mbb \in \bbz_{\geq 0} } \ \alpha^r_\mbb\, x_1^{b_1}\cdots x_k^{b_k} \) in \( k\) variables \( x_1, \ldots, x_k\), with \( \alpha^r_\mbb \in \cur{r+|\mbb|}{\bbp(\mbb)} \). By doing so, we can  package all the fundamental triples
\( \bbf_n = (\simplex{n}, \Theta_n, W_n)\)
into three elements 
\begin{equation}
\label{eq:Ftriple}
\simplex{}  := \ \sum_{b\geq 0} \simplex{b} \, t^b \ \in \ \scrI^{0}_\Pi(1|\infty),    \quad \ \ 
\Theta  := \sum_{b\geq 0} \Theta_b \, t^b \ \in \ F^0\scrN^0_\Pi (1|\infty) \quad \text{ and }\quad 
\bbw   := \sum_{b\geq 0 } W_b \, t^b \ \in \ \scrN^{-1}_\Pi(1|\infty). 
\notag
\end{equation}

In this representation, the expressions for the exterior derivatives of the currents in the fundamental triples, given in \cite[\S 3]{MR4498559},  are equivalent to saying that
\begin{equation}
\label{eq:Dtriple}
D \simplex{} = 0,\quad  D \Theta = 0, \quad \text{and} \quad D \bbw = \Theta - \simplex{}. 
\end{equation}
In other words,   we can consider the inclusions of complexes
\begin{equation}
%    \label{eq:incl_cxs}
\begin{tikzcd}
\scrI_\Pi^*(k|\infty) \ar[r, "\iota", hook] &  \scrN^*_\Pi(k|\infty) & 
F^0\scrN^*_\Pi(k|\infty) \ar[l,"\epsilon"',hook']
\end{tikzcd}
\end{equation}
and introduce the shifted cone
\begin{equation}
 %   \label{eq:uni_del_cur}
    (\bbz_\Pi^*(k|\infty), \hat D) := \mathsf{Cone}\left( \scrI_\Pi^*(k|\infty)\oplus F^0\scrN^*_\Pi(k|\infty) \xrightarrow{ \epsilon - \iota} \scrN^*_\Pi(k|\infty) \right)[-1].
\end{equation}
Then \eqref{eq:Dtriple} equivalent to the statement that 
the element \( \bbf := (\simplex{}, \Theta, \bbw)\) is a cycle in \( (\bbz^0_\Pi(1|\infty), \hat D) \) since, by definition, \( \hat D \bbf := (D\simplex{}, D\Theta, \Theta - \simplex{} - D\bbw) = (0,0,0). \)

Although arbitrary products of currents on a manifold \( M \) are not always defined, there exists an exterior product of currents \( S\times T \in \cur{r+s}{M\times N} \), for \( S \in \cur{r}{M} \) and  \( T \in \cur{s}{N} \), preserving normal and integral currents, which in turn induces an exterior product
\begin{align}
\label{eq:ext_prod}
\boxtimes \ \colon \ \unip{r}{k} \otimes \unip{s}{\ell} & \longrightarrow \unip{r+s}{k+\ell}
\end{align}
satisfying \( D(\alpha\boxtimes \beta) = (-1)^s (D\alpha)\boxtimes \beta \ + \ \alpha \boxtimes (D\beta) \), so that \( ( \unip{*}{\bullet}, D, \boxtimes) \)   becomes an associative, differential (bi)graded algebra; see Proposition \ref{prop:ext_prod}.

The  product \( \boxtimes \) induces canonical exterior products between cone complexes, as explained in \S \ref{subsec:pairings}. In broad terms, for each \( t \in \bbc\) there is a pairing of complexes
\(
\hboxt \ \colon\  \bbz_\Pi^r(k|\infty) \otimes \bbz_\Pi^s(\ell|\infty) \longrightarrow \bbz_\Pi^{r+s}(k+\ell|\infty)
\), whose definition uses the exterior product \( \boxtimes\) and certain convex combinations of currents. These pairings are homotopy-associative and their homotopy classes do not depend on \( t \in \bbc\). Subsequently, we recall in Remark \ref{rem:cup} how the exterior and cup products in Deligne cohomology are related to  \( \hboxt \).

With the preliminaries above, the quest for  multiplicative properties begins. Just like in ordinary cohomology, the cross product product  \( \alpha \times_\text{\tiny CH} \beta \) 
in higher Chow groups is  defined using the Eilenberg-Zilber maps \( \psi_{m,n} \colon \Delta^{m+n} \to \Delta^m \times \Delta^n \) in the free 
\( \mathfrak{Ab} \)-category \( \bbz \mathsf{Sch}_{\bbc}\).
It turns out that the complex \( \unip{*}{\bullet} \), and  associated subcomplexes, behaves particularly well under an extension  of the Eilenberg-Zilber maps, introduced in Definition \ref{def:psi-Pi}, thanks to the properties of current transforms. \medskip

\noindent{\bf Proposition \ref{prop:Psi12}.}{\it \ \ 
The extended Eilenberg-Zilber maps  
\[
\Psi^1 \colon \scrN^*_\Pi(1|\infty) \to\scrN^*_\Pi(2|\infty) \quad\quad \text{and} \quad  \quad 
\Psi^2 \colon \scrN^*_\Pi(2|\infty) \to\scrN^*_\Pi(3|\infty) 
\] 
are morphisms of complexes.  Furthermore, $\Psi^2\circ \Psi^1 = 0.$
}
\medskip

It is evident from its definition that  the morphism \( \Psi^1\) is related to  exterior products, as the geometric origin of the maps \( \psi_{m,n} \) was precisely to give a triangulation of the product of two topological simplices of dimensions \( m\) and \(n\), respectively. In the present context, this is can be written in terms of the exterior product \( \boxtimes \); see \eqref{eq:ext_prod}.
\medskip

\noindent{\bf Proposition \ref{prop:pf-delta}.}{\it \ \ 
 The 	simplex element \( \simplex{} \in \scrI_\Pi^0(1|\infty) \) satisfies 
 \[
 \Psi^1( \simplex{} ) \ = \ \simplex{}\boxtimes \simplex{} \ \in \scrI_\Pi^0(2|\infty) .
 \]
}
\medskip

Remarkably, the polar element \( \Theta \in F^0 \scrN_\Pi^0(1|\infty) \) has a similar behavior under \( \Psi^1\), stemming from  combinatorial relations that are explained in \S \ref{subsec:char}.  In a nutshell, we  consider the  polynomial ring 
\( \bbz[\mbx] \) on a countably infinite set of variables \( \mbx = \{ x_1, x_2, \ldots \} \)  and let  \( \bbq(\mbx) \) be its field of fractions. Then, for \( N\in \bbn\) we introduce   the rational functions
\[
{\chi}_ N(\mbx) \ =  \
\chi_N(x_1, \ldots, x_N) : = \frac{1}{x_1(x_2-x_1) \cdots (x_N - x_{N-1}) (1-x_N)} \ \in \ \bbq(\mbx).
\]
It turns out that  these elements satisfy certain \emph{shuffle relations} leading to the  result below, whose present formulation is due to Marcelo Aguiar.
\medskip

\noindent{\bf Theorem \ref{thm:charHopf}.}{\it \ \ 
The collection  \( \{ \chi_N(\mbx) \mid N \in \bbn\} \) defines a character \( \bchi \) of the permutation Hopf algebra \( \perHA\) with values in the field \( \bbq(\mbx)\), 
defined by the assignment
\[
\bchi \ \colon \ \scrF_\sigma \ \longmapsto \chi_m\left(x_{\sigma(1)}, \ldots, x_{\sigma(m)}\right), \quad \sigma \in \symg{m},
\]
where \( \{  \scrF_\sigma \mid \sigma \in \symg{m}, m\geq 0 \} \) is the \( \scrF \)-basis for \( \perHA\).
}
\medskip

Now, we can determine the behavior of the  next  element in the fundamental triple.
\medskip

\noindent{\bf  Theorem \ref{thm:pf-theta}.}{\it  \ \ 
The 	polar element \( \Theta \in F^0\scrN_\Pi^0(1|\infty) \) satisfies 
 \[
 \Psi^1( \Theta ) \ = \ \Theta\boxtimes \Theta \ \in F^0 \scrN_\Pi^0(2|\infty) .
 \]
}

From this theorem, we conclude in Corollary \ref{cor:cycHCG} that 
\[ \Psi^1 (\bbf) := ( \Psi^1 \simplex{},\  \Psi^1 \Theta,\  \Psi^1 \bbw) = (\simplex{}\boxtimes \simplex{},\  \Theta\boxtimes \Theta,\  \Psi^1 \bbw)  \ \in \ \bbz^0_\Pi(2 | \infty ) 
\]
 is a cycle in the complex \( \bbz^*_\Pi(2 |\infty) \).

In order to establish the relation between products in the cohomology theories under study  we need to compare the  cycles   \( \Psi_1(\bbf)\  \) and \(\ \bbf\hboxt\bbf \) in \( \bbz_\Pi^*(2 | \infty)\). This comparison follows from the calculation of the homology of \( \scrN_\Pi^*(2|\, \infty) \), which goes as follows. 
\medskip 

\noindent{\bf Theorem \ref{thm:coh-2}.}{\it \ \
The  product complex $\scrN^*_\Pi(2|\infty)$ is quasi-isomorphic to $\underline{\bbc}$ (concentrated at $0$). In particular,
$H^{-1}(\scrN^*_\Pi(2|\infty)= 0.$
}
\medskip 
 
\noindent{\bf Corollary \ref{cor:homology}.}{\it \ \ 
The elements 
\( 
\Psi^1(\bbf ) = (\simplex{}\boxtimes \simplex{}, \Theta \boxtimes \Theta, \Psi^1(\bbw) ) \)\ and \ \( \bbf \hboxt \bbf  = (\simplex{}\boxtimes \simplex{}, \Theta \boxtimes \Theta, \msfP_t(\bbf, \bbf) ) 
\)
represent the same homology class in \(  H^0 ( \bbz_\Pi^*(2 | \infty)) \). 
}
\medskip

Once this corollary is established, the proofs of the main results\textemdash specifically Theorem \ref{thm:MAIN} and Corollary \ref{cor:ext}\textemdash outlined in the beginning of the introduction, become a tautology due to the inherent properties of the current transforms.
\medskip

\noindent{\sl Acknowledgements:}\ \ The author would like to thank Marcelo Aguiar for his invaluable insight, which was instrumental for the elegant formulation and proofs in \S \ref{subsec:char}.

%%%%%%%%%%% \input{Multicomplexes.tex}
% !TEX root =  main_HCG_prod.tex

\section{Multi-graded complexes of currents and the fundamental triple}
\label{sec:multi}

In this preliminary section we repackage, in a succinct form,  the complexes of currents that appear in the construction of the regulator maps in \cite{MR4498559} and recall the definition and properties of the fundamental triple of currents \( \bbf_n = (\simplex{n}, \Theta_n, W_n) \).

The basic objects are  \emph{multi-graded} \( (k+1)\)-complexes of modules over a commutative ring with commuting differentials, along with their associated total product complexes. They all have the form
\[
\Gamma^{*,\bullet}= \left(
\Gamma^{a,\mbb},\, d, \, \delta^j \ \mid \ a \in \bbz,\  \mbb \in \bbz^k_{\geq 0},\,  j = 1, \ldots, k
\right),
\]
with \(  d \colon \Gamma^{a,\mbb} \to \Gamma^{a+1, \mbb} \)\quad
and\quad 
\(  \delta^j \colon \Gamma^{a,\mbb}(k) \to \Gamma^{a,\mbb+ \mbe_j}(k) \), 
where \(  \mbe_j = (0, \ldots, 1,\ldots, 0) \) is the \(  j \)-th canonical basis element. In summary,
\(\ 
d^2=0,\quad (\delta^j)^2 =0, \quad  d \circ \delta^j = \delta^j \circ d, \quad \delta^i\circ \delta^j = \delta^j \circ \delta^i,\quad 
\text{ for all }  1\leq i, j \leq k.
\)
\begin{definition}
\label{def:tot}
Given a \(  (k+1) \)-complex \(  (\Gamma^{*,\bullet}, d, \delta^j)  \) as above,   denote
\begin{equation}
\label{eq:tot}
\Gamma^r_\Pi := \prod_{a+|\mbb|=r}\Gamma^{a, \mbb},  \quad\ \    r \in \bbz,
\end{equation}
where \(  |\mbb| = b_1 + \cdots + b_k \) for  a \(  k \)-tuple   \(  \mbb = (b_1, \ldots, b_k)  \in \bbz_{\geq 0}^{ k}  \).
It is convenient to write  \(  \alpha \in \Gamma^r_\Pi \) as a formal power series 
\(
\alpha = \sum_{\mbb  } \alpha^r_\mbb\ x_1^{b_1} \cdots x_k^{b_k},
\)
%on \( k\)-variables \(  x_1, \ldots, x_k \), 
with \( \alpha_\mbb^r \in \Gamma^{r-|\mbb|, \mbb}. \) 
\ Consider
\begin{equation}
\label{eq:dalpha}
d\alpha \ := \ \sum_{\mbb} \left( d \alpha_\mbb^r \right) \ x_1^{b_1} \cdots x_k^{b_k} \ \in \ \Gamma^{r+1}_\Pi
\quad \text{ and } \quad 
\hat \delta \alpha \ : = \ \sum_\mbb \left( \hat \delta \alpha\right)^{r+1}_\mbb \ x_1^{b_1} \cdots x_r^{b_k} \ \in \ \Gamma^{r+1}_\Pi,
\end{equation}
where 
\(
\left( \hat \delta \alpha \right)^{r+1}_\mbb \ := \ 
\sum_{j=1}^k (-1)^{b_j+\cdots + b_k} \delta^j\left( \alpha^r_{\mbb-\mbe_j} \right).
\)
It is easy to see that \(  d^2 = 0, \ \hat \delta^2 = 0 \) and
\(  \hat \delta \circ d = d \circ \hat \delta. \ \)  Now, one defines the total differential
\begin{equation}
\label{eq:D}
D \ \colon\ \Gamma^r_\Pi \longrightarrow \Gamma^{r+1}_\Pi
\end{equation}
as \(  D \alpha = d\alpha + (-1)^r \hat\delta \alpha,\)
and obtains the \emph{total product complex} 
\(  \left( \Gamma^*_\Pi, D \right). \)
\end{definition}
\subsection{Universal multi-graded complexes of currents}
\label{subsec:univ}

Given a category \(  \scrC  \),  let   \(  \bbz\scrC \) denotes the free \(  \mathsf{Ab} \)-category  defined by
\(
%\begin{equation}
%\label{eq:ZC}
\mathsf{Ob}(\bbz\scrC) = \mathsf{Ob}(\scrC)  \  \text{ and } \   \Hom{\bbz\scrC}{X}{Y} := \zhom{\scrC}{X}{Y},
%\end{equation}
\)
where \( \zhom{\scrC}{X}{Y} \) is the free abelian group on \(  \Hom{\scrC}{X}{Y}. \)
%When working in  the category \(  \sch{\bbc} \) of schemes of finite type over \( \bbc \) and   no confusion is likely to arise, we use  \(  X \) to denote both a scheme \(  X \) and the complex analytic variety \( X(\bbc)^\text{an}   \) given by \(  \bbc \)-valued points of \(  X \) with the analytic topology. 

Let \( [\mbz]=[z_0:\cdots : z_b] \) be the homogeneous coordinates in \( \bbp^b, \ b\in \bbn,  \) and let 
\(  \iota_r \colon \bbp^{b-1} \hookrightarrow \bbp^b \) be the inclusion of \(  \bbp^{b-1} \) in \(  \bbp^b \) as the \(  r\)-th coordinate hyperplane, given   by \(  z_r = 0. \)  Also, let \(  H_\infty(b) \subset \bbp^b \) be the hyperplane at infinity defined by \(  z_0 + \cdots + z_b = 0 \) and let \(  \Delta^b = \bbp^b - H_\infty(b) \) be the algebraic \(  b \)-simplex.
The alternating sum of the \emph{face maps} 
\( \iota_r\) gives the \emph{boundary element} \( \delta := \sum_{r=0}^b (-1)^r \iota_r\ \in \  \zhom{\sch{\bbc}}{\bbp^{b-1}}{\bbp^b}. \) 

\begin{notation}
\label{not:spaces}
Given a \(  k \)-tuple \(  \mbb \in \bbz^{  k}_{\geq 0}\), let
\(\bbp(\mbb) := \bbp^{b_1} \times \cdots \times \bbp^{b_k}
\)
be the corresponding product of  projective spaces.  Set \( \bbp(\mbb) = \emptyset\) for \( \mbb \in \bbz^k - \bbz^k_{\geq 0} \), 
and denote 
\begin{equation}
\label{eq:gen_delta}
\delta^i := (1 \times \cdots \times \delta \times \cdots \times 1) \in \zhom{\sch{\bbc}}{\bbp(\mbb - \mbe_i)}{\bbp(\mbb)}. 
\end{equation}

\end{notation}

\begin{definition}
\label{def:dbl_cx}
Fix \(  k\geq 0 \). Given $a \in \bbz$ and $\mbb=(b_1, \ldots, b_k) \in \bbz_{\geq 0}^k $, define
\begin{equation}
\uni{a,\mbb}(k) := \cur{2|\mbb|+a}{\bbp(\mbb)} = 
\curd{-a}{\bbp (\mbb)},
\end{equation}
as the group of de Rham currents of dimension \(  -a \geq 0 \) in \(  \bbp(\mbb), \) with $\uni{a,\mbb}(k) = 0$ whenever   $a>0$. The collection $\{ \uni{a,\mbb}(k) \mid a \in \bbz, \mbb \in \bbz_{\geq 0}^k \}$ becomes a \(  (k+1) \)-complex as follows. 
The   exterior derivative of currents gives  $d \colon \uni{a,\mbb}(k) \to \uni{a+1, \mbb}(k) $, and we define $ \delta_\#^i \colon \uni{a,\mbb - \mbe_i}(k) \to \uni{a,\mbb}(k)$ as the map of currents induced by the morphism  \(  \delta^i  \) introduced in   Notation \ref{not:spaces}.
It is clear that
\[
d^2 = 0,\quad  (\delta^i)^2 = 0,\quad   d\circ \delta_\#^i = \delta^i_\# \circ d, \quad \text{ and }\quad    \delta^i_\# \circ \delta^j_\# = \delta^j_\#\circ \delta^i_\#,\quad  \text{ for all } \   0\leq i, j \leq k. 
\] 
Hence
\( ( \uni{*,\bullet}(k); d,\ \delta_\#^1, \ldots, \delta_\#^k ) \) is a \(  (k+1) \)-complex as in Definition \ref{def:tot}, 
and  its associated {total product complex}  \(  (\unip{*}{k}, D) \) is defined as in  \eqref{eq:tot} and \eqref{eq:D}.
\end{definition}

\subsection{The exterior product}

Using the exterior product of currents \cite[\S 4.1.8]{MR0257325} we introduce pairings between the  complexes \(   \unip{*}{\bullet}\)  introduced above. 

\begin{definition}
\label{def:exterior}
Given $\balpha \in \unip{r}{k}\) and  \(\bbeta \in \unip{s}{\ell},$ define
\[
\balpha \boxtimes \bbeta := 
\sum_{\mbd=\mba*\mbb\, \in\, \bbz_{\geq 0}^{k+\ell}}  (-1)^{r|\mbb|}  \left(  \alpha^r_\mba \times \beta^s_{\mbb}\right)\, t_1^{d_1}\cdots t_{k+\ell}^{d_{k+\ell}} ,
\]
where $\mba*\mbb = (a_1, \ldots, a_k, b_1, \ldots, b_\ell)$ is the concatenation of indices and $\alpha^r_\mba \times \beta^s_\mbb \in \cur{r+s+|\mba|+|\mbb|}{\bbp(\mba*\mbb)}$
is the exterior product of currents.
\end{definition}

\begin{proposition}
\label{prop:ext_prod}
The pairing
\(
\boxtimes \ \colon \ \unip{*}{k}\otimes \unip{*}{\ell} \longrightarrow \unip{*}{k+\ell}
\)
is associative and satisfies
\[
D(\balpha \boxtimes \bbeta) \ = \ (-1)^s (D \balpha) \boxtimes \bbeta \ + \ \alpha \boxtimes (D \bbeta),
\]
when $\balpha \in \unip{r}{k}$ and $\bbeta \in \unip{s}{\ell}$. Therefore the  complex $\left( \unip{*}{\bullet}:= \bigoplus_k \unip{*}{k}, \ D, \ \boxtimes \right)$ becomes  a differential (bi)graded algebra under the product \( \boxtimes\).

\end{proposition}
\begin{proof}
This follows directly from \cite[\S 4.1.8]{MR0257325} and the definition of \( D\). 
\end{proof}

\subsection{Geometric subcomplexes}\label{subsec:interest}  
\label{subsec:subcxs}

There are particular classes of currents that behave well under the transform operations, which are also used to describe Deligne-Beilinson cohomology in \cite{MR4498559}. Here  we introduce the corresponding \(  (k+1) \)-subcomplexes of $ \uni{*,\bullet}(k)$.   
 
Consider a complex manifold $M$ and let $H\subset M$ be  a complex submanifold. A  \emph{locally normal} current $T \in \scrN_\text{loc}^{k}(M)$   is said to \emph{vanish suitably along $H$} if:
\begin{enumerate}[1)]
    \item $\norm{T}(H)=0$, where \( \norm{T}\) is the Radon measure associated to \( T\);
    \item $\norm{dT}(H)=0$, where \( \norm{dT}\) is the Radon measure associated to \( dT\);
    \item For a fixed Riemannian metric on $M$, there is some $\epsilon >0$ such that for each tubular neighborhood $W^\tau$ of radius $\tau<\epsilon$ and smooth oriented boundary $\partial W^\tau$, the intersection $T \cap \llb \partial W^\tau \rrb$ exists.
\end{enumerate}

\begin{example}
It is shown in \cite{MR4498559} that the currents $\Theta_n,\,  \simplex{n}$ and $ W_n$ (see Definition \ref{def:triple}) are normal and vanish suitably along the hyperplane $H_\infty(n) \subset \bbp^n.$
\end{example}

\begin{definition}
\label{def:vanish_infty}
Given $\mbb = (b_1, \ldots, b_k) \in \bbz_{\geq 0}^k$, consider the divisor with simple normal crossings  
\( 
H_\infty(\mbb) = \bigcup_{j=1}^k \left( \bbp^{b_1}\times \cdots \times H_\infty(b_j) \times \cdots \times \bbp^{b_k}\right) \subset \bbp(\mbb)
\)
and denote by $\nor{k}{\bbp(\mbb) | \infty}$ the group of normal currents that vanish suitably along each smooth component of the divisor. These currents form a complex
$\left( \scrN^*(\bbp(\mbb)|\infty), d\right)$ under the usual exterior differentiation of currents. Furthermore,  the differentials 
$\delta_j \colon \cur{k}{\bbp(\mbb)} \to \cur{k+2}{\bbp(\mbb + \mbe_j)}$, introduced above,    preserve the subgroups $\scrN^*(\bbp(\mbb)|\infty)$.
In fact, the subgroups below are all closed under these differentials.
\[
\begin{alignedat}{3}
&    \msfF^p\mathbb{\Gamma }^{a, \mbb}(k) && := F^{p+a+|\mbb|}\  
   \curd{a+|\mbb|}{\bbp(\mbb)} &&,   \text{ where } \ F^q \cur{*}{-} \text{\  is the Hodge filtration.} \\
&    \scrN^{a,\mbb}(k) && =  \scrN^{2|\mbb|+a}(\bbp(\mbb) ) && =  \text{ normal curr. of dimension} -a \text{ on } \bbp(\mbb). \\
&  \scrN^{a,\mbb}(k| \infty) && = \scrN^{2|\mbb|+a}(\bbp(\mbb)| \infty ) && =  \text{normal curr.  \emph{vanishing suitably} at } H_\infty(\mbb).  \\
&    \scrI^{a,\mbb}(k) && =  \scrI^{2|\mbb|+a}(\bbp(\mbb) ) && =  \text{ integral curr. of dimension} -a \text{ on } \bbp(\mbb).   \\
&    \scrI^{a,\mbb}(k| \infty) && = \scrI^{a,\mbb}(k) \cap \scrN^{a,\mbb}(k| \infty) && =  \text{integral curr.  \emph{vanishing suitably} at } H_\infty(\mbb).
\end{alignedat}
\]
It follows that  the  total product complex
$\left( \scrN_\Pi^*(k|\infty), D\right)$  is a subcomplex of $(\unip{*}{k}, D).$
Similarly, the subgroups listed above yield   corresponding subcomplexes of
\( \unip{*}{k}\):
\[
\begin{tikzcd}[column sep=0.8cm]
&  \scrI^{*}_\Pi(k| \infty) \ar[r, hook] \ar[d, hook]  &  \scrI^{*}_\Pi(k) \ar[d, hook]    \\
\cdots \subset  \msfF^{-1}\scrN^{*}_\Pi(k| \infty) \subset \msfF^{0}\scrN^{*}_\Pi(k| \infty) \subset \cdots \subset \msfF^{p}\scrN^{*}_\Pi(k| \infty) \subset  \cdots  
\subset \hspace{-1cm} & \scrN^{*}_\Pi(k| \infty) \ar[r, hook]  &  \scrN^{*}_\Pi(k)  .
\end{tikzcd}
\]
\end{definition}

\subsection{Cone complexes and multiplication pairings} 
\label{subsec:pairings}
Using the inclusions 
\[
\begin{tikzcd}
\scrI_\Pi^*(k | \infty) \ar[r, hook, "\ve"]
 &
 \scrN_\Pi^*(k | \infty)
 &
 \msfF^0 \scrN_\Pi^*(k | \infty) \ar[l, hook', "\iota"'],
 \end{tikzcd}
 \]
for each \( k\geq 0 \)   define  the complex
\begin{equation}
    \label{eq:uni_del_infty}
    \bbz_{\Pi}^*(k|\infty) = \mathsf{Cone}\left( \scrI_\Pi^*(k|\infty)\oplus \msfF^0
    \scrN^{*}_\Pi(k|\infty)\xrightarrow{\ \  \epsilon - \iota\ \ } \scrN^*_\Pi(k|\infty) \right)[-1].
\end{equation}
For simplicity,   write \( \hat D(A) = (D\balpha_1, \ D\balpha_2, \ \balpha_2 -  \balpha_1 - D \balpha_3)\) for the differential of an element \( A=(\balpha_1, \balpha_2, \balpha_3)\in   \bbz_{\Pi}^r(k|\infty) \), omitting the inclusions \( \ve\) and \( \iota\).
 
 \begin{definition}
 \label{def:cup}
 Fix \( t \in \bbc\).  Given  \( \balpha, \balpha'  \in \scrN^r_\Pi(k) \), denote
\(
 \Xi_t(\balpha, \balpha') := t \balpha \ + \ (1-t) \balpha'.
 \)
Then, for \( A = (\balpha_1, \balpha_2, \balpha_3) \in \bbz_\Pi^r(k|\infty) \) and \( B = (\bbeta_1, \bbeta_2, \bbeta_3) \in \bbz_\Pi^s(\ell|\infty) \), define
\begin{equation}
\label{eq:Pt}
\msfP_t(A,B) := (-1)^s \balpha_3 \boxtimes \Xi(\bbeta_1, \bbeta_2)\ + \ \Xi(\balpha_1, \balpha_2)\boxtimes \bbeta_3  \ \ \in \ \scrN^{r+s -1}_\Pi(k+\ell).
\end{equation}
Now, define a pairing 
 \begin{align}
 \label{eq:cupt}
\hboxt \ \colon\  \bbz_\Pi^r(k|\infty) \otimes \bbz_\Pi^s(\ell|\infty) & \longrightarrow \bbz_\Pi^{r+s}(k+\ell|\infty)\\
 A  \otimes B \hspace{0.64cm} & \longmapsto (\balpha_1\boxtimes \bbeta_1,\ \balpha_2\boxtimes \bbeta_2,\ \msfP_t(A,B)). \notag
 \end{align}
 \end{definition}

\begin{remark}
\label{rem:dbc}
The pairings \( \hboxt \)   above are directly related to the product in Deligne-Beilinson cohomology. See Remark \ref{rem:cup}.

\end{remark}

 \begin{proposition}
 \label{prop:cupt}
 For each \( t \in \bbc\),  the pairing \(  \hboxt \) satisfies the following properties.
 \begin{enumerate}[i.]
 \item If \( A \in \bbz_\Pi^r(k|\infty) \) and \( B \in \bbz_\Pi^s(\ell|\infty) \) then
 \[ \hat D (A\ \hboxt B) = (-1)^s ( \hat D A) \hboxt  B \ + \ A\, \hboxt (\hat D B).\]
 It follows that  \( \hboxt \colon \bbz_\Pi^*(k|\infty) \otimes \bbz_\Pi^*(\ell|\infty) \to 
 \bbz_\Pi^*(k+\ell|\infty) \) defines a pairing of complexes.
 \item The maps \( \hboxt \) are   homotopic to each other, for all \( t \in \bbc\).
 \item \( \hboxt \) is homotopy-associative.
\end{enumerate} 
 \end{proposition}
 \begin{proof}
 The first assertion from directly from the claim below. 
 
\noindent{\bf Claim.}
Given \( A = (\balpha_1, \balpha_2, \balpha_3) \in \bbz^r_\Pi(k | \infty ) \) and \( B = (\bbeta_1, \bbeta_2, \bbeta_3) \in \bbz_\Pi^s(\ell  | \infty ) \) one has:
\[ D\msfP_t(A, B) \ +\  (-1)^s \msfP_t(\hat DA, B)\  +\  \msfP_t(A, \hat D B)  \ = \ 
\balpha_2\boxtimes \bbeta_2 \ - \ \balpha_1\boxtimes \bbeta_1. \]
\begin{proof}[Proof of Claim]
By definition, the three terms on the left hand side of the identities are the following:
\begin{align*}
D\msfP_t (\balpha, \bbeta) & =
 D \left\{  (-1)^s \balpha_3 \boxtimes \Xi_t (\bbeta_1,\bbeta_2) \ +\ \Xi_{(1-t)}(\balpha_1,\balpha_2) \boxtimes \bbeta_3   \right\} \\
 & =
 (-1)^s\{ (-1)^s (D\balpha_3)\boxtimes \Xi_t(\bbeta_1,\bbeta_2) + \balpha_3\boxtimes \Xi_{(1-t)}(D\bbeta_1, D\bbeta_2) \} \\
 &  +  (-1)^{s-1} \Xi_{(1-t)}(D\balpha_1, D\balpha_2)\boxtimes \bbeta_3 \ + \ \Xi_{(1-t)}(\balpha_1,\balpha_2)\boxtimes D\bbeta_3\\
 & =
 (D\balpha_3)\boxtimes \Xi_t(\bbeta_1,\bbeta_2) \ + \ (-1)^s \balpha_3\boxtimes \Xi_{t}(D\bbeta_1, D\bbeta_2)  \\
 &  +  (-1)^{s-1} \Xi_{(1-t)}(D\balpha_1, D\balpha_2)\boxtimes \bbeta_3 \  + \ \Xi_{(1-t)}(\balpha_1,\balpha_2)\boxtimes D\bbeta_3, 
\end{align*}
and
\begin{align*}
 (-1)^s \msfP_t  (\hat D \balpha, \bbeta)  &= (-1)^s\left\{ (-1)^s ( \balpha_2 - \balpha_1 - D\balpha_3) \boxtimes \Xi_t(\bbeta_1,\bbeta_2)  +  \Xi_{(1-t)}(D\balpha_1,D\balpha_2)\boxtimes \bbeta_3
\right\} \\
& =
( \balpha_2 - \balpha_1) \boxtimes \Xi_t(\bbeta_1,\bbeta_2)  -  ( D\balpha_3) \boxtimes \Xi_t(\bbeta_1,\bbeta_2) + (-1)^s \Xi_{(1-t)}(D\balpha_1,D\balpha_2)\boxtimes \bbeta_3,
\end{align*}
and
\begin{align*}
 \msfP_t  (\balpha, \hat D\bbeta) & = (-1)^{s+1} \balpha_3\boxtimes \Xi_{t}(D\bbeta_1,D\bbeta_2) + \Xi_{(1-t)}(\balpha_1, \balpha_2) \boxtimes (\bbeta_2 - \bbeta_1 - D\bbeta_3) \\
& = (-1)^{s+1} \balpha_3\boxtimes \Xi_{t}(D\bbeta_1,D\bbeta_2) +  
\Xi_{(1-t)}(\balpha_1, \balpha_2) \boxtimes (\bbeta_2 - \bbeta_1)  -  \Xi_{(1-t)}(\balpha_1, \balpha_2) \boxtimes (D\bbeta_3).
\end{align*}
Adding the identities above proves the claim:
\begin{align*}
D\msfP_t(\balpha, \bbeta) \ & +\  (-1)^s \msfP_t(\hat D\balpha, \bbeta)\  +\  \msfP_t(\balpha, \hat D \bbeta)   = 
( \balpha_2 - \balpha_1) \boxtimes \Xi_t(\bbeta_1,\bbeta_2) + \Xi_{(1-t)}(\balpha_1, \balpha_2) \boxtimes (\bbeta_2 - \bbeta_1)  \\
& =
( \balpha_2 - \balpha_1) \boxtimes \{t\bbeta_1 + (1-t)\bbeta_2\} +  \{ (1-t)\balpha_1 + t \balpha_2\}  \boxtimes (\bbeta_2 - \bbeta_1)   = \balpha_2\boxtimes \bbeta_2  \ - \  \balpha_1 \boxtimes \bbeta_1.
\end{align*}
\end{proof}
 
The second statement in the proposition follows from a standard argument. As to the last statement, since it is not used throughout the paper, we leave it as an exercise for the reader.
 \end{proof}

\begin{example}
\label{exmp:cycDel}
For each $t\in \bbc$,  the element
 \(\, \bbf\hboxt \bbf = (\simplex{} \boxtimes \simplex{}, \ \Theta\boxtimes \Theta, \  \msfP_t (\bbf, \bbf))\) is a cycle in the complex \( \bbz_\Pi^0(2|\infty) \), where 
\( \bbf = (\simplex{}, \Theta, \bbw) \) is the fundamental triple introduced in next section. 
In particular, 
\(
%\begin{equation}
%\label{eq:Pt}
  D \msfP_t (\bbf, \bbf) \ = \ \Theta\boxtimes \Theta \ - \ \simplex{} \boxtimes \simplex{}.
%\end{equation}
\)
 \end{example}

\begin{remark}
It follows from Corollary \ref{cor:cycHCG} and Example \ref{exmp:cycDel}   that  both 
\( (\simplex{} \boxtimes \simplex{}, \ \Theta\boxtimes \Theta, \  \msfP_t (\bbf, \bbf)) \)
and 
\(  (\simplex{}\boxtimes \simplex{}, \Theta\boxtimes \Theta, \Psi^1 \bbw) \) are cycles in \(\bbz^0_\Pi(2 | \infty )\). 
We show in the next section that these cycles are homologous, and the proof of the main result is  a  consequence of this fact. 
\end{remark}

\subsection{Revisiting the Fundamental Triple of Currents}
\label{subsec:FTC}

Here we recall the construction of the triple \( \bbf_n = (\simplex{n}, \Theta_n, W_n) \) introduced in \cite{MR4498559}, and summarize its properties.

First, consider the semi-algebraic set \( S_j  \subset \bbp^n \),  \(  1\leq j \leq n \),  defined as
\(
 S_j := \{ [z_0:\cdots:z_n] \ \mid \ z_j = t ( z_0 + \cdots + z_j) , \text{ for some } t \in [0, 1] \}.
\)
%with Definition \ref{def:sacurr},
For each \(  0\leq j < n \) it is shown that one has a proper intersectioan 
\( R_{n,j}  = S_{n}  \cap S_{n-1} \cap \cdots \cap   S_{j+1} \subset \bbp^n. \) Let   \( \llb R_{n,j} \rrb\) denotes the associated semi-algebraic chain,  oriented as in \cite[\S3.1]{MR4498559}.
Now, set \(  \tet{0}{n}  =1 \) and for \(0< j \leq n \) define meromorphic \(  j \)-forms   \(  \tet{j}{n}  \) in \(  \bbp^n \) by  
   \begin{equation}
   \label{eq:theta}
 \tet{j}{n}  := \sum_{r=0}^j (-1)^r \frac{dz_0}{z_0}\wedge \cdots \wedge \widehat{\frac{dz_r}{z_r}} \wedge \cdots \wedge \frac{dz_j}{z_j}.
 \end{equation}
%This form has log-poles along the divisor   \om{0}{n}   = 0\), and  \( \om{j}{n}  := (-1)^j\hbar_j\,  \tet{j-1}{n} \ \)
for \(  1 \leq j \leq n \),
where
\[    
\hbar_j[\mbz] =\begin{cases}
\vspace{0.1in}
 \log \left( 1 - \frac{\ve_j(\mbz)}{z_j}\right)  
&, \text{ if } [\mbz] \notin S_j \\
\vspace{0.1in}
0 &, \text{ if }  [\mbz] \in S_j.
\end{cases} 
\]
Here \(  \ve_j(\mbz) = z_0 + \cdots + z_j\)\,   and \( \log ( - ) \) denotes the principal branch of the logarithm.

\begin{definition}
\label{def:triple}
The \emph{fundamental triple} \( \bbf_n = (\simplex{n}, \Theta_n , W_n)\)   of currents on  \( \bbp^n\) is defined as follows. 
\begin{enumerate}[a.]
\item Let $\simplex{n}\subset \Delta^n \subset \bbp^n$ be the topological simplex. For simplicity,    the symbol \( \simplex{n}\) is also used to denote the corresponding integral current $\llb \simplex{n} \rrb \in \intcur{n}{\bbp^n}\subset \mathbb{\Gamma}^{-n,n}(1),$ with its usual orientation.
\item Define
\begin{equation}
    \label{eq:hatTheta}
    \Theta_n := \frac{ (-1)^{\frac{n(n+1)}{2}}}{\twopii^n} [\theta_n] \ \in F^n\scrN^{n}(\bbp^n) \subset \msfF^0\mathbb{\Gamma}^{-n,n}(1),
\end{equation}
where $[\theta^n]$ is the current represented by  $\theta_n $, that is,  $ [\theta_n]  \colon \varphi \mapsto \int_{\bbp^n} \theta_n \wedge \varphi .$
\item Set
\(  W_0  = 0 \) and for \(n \geq 1  \) define
\begin{equation}
\label{eq:Wn}
W_n   := (-1)^{n} \sum_{j = 1}^{n} \frac{(-1)^{\binom{j}{2}}}{\twopii^{j}}\  \llb R_{n,j}  \rrb\wedge \om{j}{n} \ \in \scrN^{n-1}(\bbp^n) \subset \mathbb{\Gamma}^{1-n,n}(1).
\end{equation}
 \end{enumerate}
 It is shown in \cite[Lem. 4.2]{MR4498559} that all three currents \( \simplex{n}, \Theta_n\) and \( W_n\) \textbf{vanish suitably at \( \infty\).}
\end{definition}

\begin{remark}
\label{rem:const}
The original definition of \( \Theta_n \) and \( W_n\) differs from the one introduced here by the factor  \( \frac{(-1)^{{n(n+1)}/{2}}}{\twopii^n}\). Under the present renormalized version,  the formulas shown in \cite[\S 3]{MR4498559}  for the exterior derivatives of these currents can be rewritten as:
\[
d\simplex{n}   +  (-1)^n \delta_{\#} \simplex{n-1} = 0, \quad  d\Theta_n  +  (-1)^n \delta_{\#} \Theta_{n-1}  = 0, \quad \text{ and }\quad 
  d W_n +  (-1)^{n-1} \delta_{\#} W_{n-1}  = \Theta_n - \simplex{n}.
\]
 See Notation \ref{not:spaces} for the definition of \( \delta_\#\).
\end{remark}

\begin{definition}
\label{defprop:triple}
Using the triples \( \bbf_n\),   define the following elements:
%\[
%\begin{alignedat}{3}
%&    \simplex{} && := \sum_{n\geq 0} \simplex{n}\ x^n &&  \in \  \scrI_\Pi^{0}(1 | \infty), \\
%&    \Theta &&  := \sum_{n\geq 0} \Theta_n\ x^n  &&  \in \ F^0 \scrN_\Pi^{0}(1 | \infty),  \\
%&    \bbw   && := \sum_{n\geq 0} W_n\ x^n \  && \in \   \scrN_\Pi^{-1}(1 | \infty),
%\end{alignedat}
%\]
\[
\simplex{}  := \sum_{n\geq 0} \simplex{n}\ x^n   \in \  \scrI_\Pi^{0}(1 | \infty), \quad\ 
    \Theta   := \sum_{n\geq 0} \Theta_n\ x^n    \in \ \msfF^0 \scrN_\Pi^{0}(1 | \infty),  \quad\ \text{and}\ \quad \   \bbw    := \sum_{n\geq 0} W_n\ x^n \   \in \   \scrN_\Pi^{-1}(1 | \infty),
\]
and call \(\ \ \bbf = ( \simplex{}, \Theta, \bbw) \ \in\  \bbz^0_\Pi(1 | \infty)  \)  the \ \emph{fundamental triple}.
\end{definition}

 \begin{proposition}
 \label{cor:cycle}
 The triple \( \bbf = (\simplex{}, \Theta, \bbw) \in \bbz^0_\Pi(1 | \infty)  \) is a cycle in the complex \( \left( \bbz^*_\Pi(1 | \infty), \hat D\right) \); see \ref{eq:uni_del_infty}. Equivalently,
\(
D\simplex{} \ = \ D \Theta \ = \ 0 \quad \text{ and } \quad D\bbw \, = \, \Theta  - \simplex{}.
\)
 \end{proposition}
 \begin{proof}
 Follows  directly from Remark \ref{rem:const},  \cite[Prop. 3.1, Prop. 3.3, Cor 3.5]{MR4498559} and the definitions.
 \end{proof}

%%%%%%%%%%%%%%%%%.  \input{EZproducts.tex}
 % !TEX root =  main_HCG_prod.tex

\section{Extended Eilberg-Zilber maps and products}

In this section we recast algebraic and topological incarnations of the classical Eilenberg-Zilber (EZ) constructions to define maps between the various product complexes of currents. Then we examine the behavior of the triple \( \bbf =( \simplex{}, \Theta, \bbw) \) under the resulting \emph{extended EZ map}.

\subsection{The Eilenberg-Zilber morphisms}

Let us first recall  the definition of the algebraic-geometric  simplicial Eilenberg-Zilber (EZ)  morphism and its ``cubical'' counterpart.
Let \(  \acube{N}:= \bba^1 \times \cdots \times \bba^1 \)  be the algebraic \(  N \)-cube, with coordinates   \(  \mbt = (t_1, \ldots, t_N) \), and let \(  \mbz = (z_0, \ldots, z_N)\) be the coordinates in   \(  \asimplex{N} \).
Using  \(  \ve_j(\mbz) = z_0 + \cdots + z_j \),  \  define  isomorphisms 
\begin{equation}
\label{eq:phi}
\phi_N \colon \asimplex{N} \xrightarrow{\  \cong\  }  \acube{N} \quad \text{and} \quad \dif{N}= \phi_N^{-1} \colon \acube{N} \xrightarrow{\  \cong\  } \asimplex{N}
\end{equation}
  by 
\(
\phi_N ( \mbz) = (\ve_0(\mbz), \ve_1(\mbz), \ldots, \ve_{N-1}(\mbz) )\) and 
\(\dif{N} (\mbt) =  (t_1, t_2-t_1, \ldots, t_N - t_{N-1}, 1 - t_N),
\)
respectively.
 
Now, consider the set  \(  \shuf{m}{n} \subset \symg{m+n}\)   of shuffles of type \(  (m, n) \), defined as the permutations \(  \tau  \)   of \(  [m+n] = \{ 1, \ldots, m+n \} \) whose restrictions to \(  [m] \) and to \(  [m+n] - [m] \) are increasing. 
Each \( \tau \in \shuf{m}{n}\) induces an isomorphism
\( \lambda_\tau\colon \asimplex{m+n} \xrightarrow{ \ \cong \ }    \asimplex{m}\times \asimplex{n} \) defined as follows. If \(
s_k \colon \asimplex{N} \to \asimplex{N-1}, \ \ k = 0, \ldots, N-1,
\)
denotes the  (co)degeneracy map  given by   \(  s_k(\mbz) = (z_0, \ldots, z_{k-1}, z_k+z_{k+1}, z_{k+2} ,\ldots, z_N)\), then
\begin{equation}
\label{eq:lambda}
\lambda_\tau (\mbz) := \left( s_{\sigma(m+1)}\circ s_{\sigma(m+2)} \circ \cdots \circ s_{\sigma(m+n)}\right) \times 
\left( s_{\sigma(1)}\circ s_{\sigma(2)} \circ \cdots \circ s_{\sigma(m)}\right).
\end{equation}
Define the  \emph{Eilenberg-Zilber (EZ)  morphism} \(\psi^\Delta_{m,n}  \in \zhom{\mathsf{Sch}}{\asimplex{m+n}}{\asimplex{m}\times  \asimplex{n}}
 \) by
  \begin{equation}
\label{eq:Psi_Delta}
\psi^{\Delta}_{m,n}\, \,  = \sum_{\tau \in \shuf{m}{n}} \epsilon(\tau)\,  \lambda_\tau,
\end{equation}
where \( \epsilon(\tau) \) is the sign of  \( \tau\).
Similarly,  using the usual action of \(  \symg{m+n} \) on \(\   \acube{m+n}  \),  define the \emph{cubical  EZ morphism}
\(  \psi^\oblong_{m,n} \in \zhom{\mathsf{Sch}}{\acube{m+n}}{\acube{m}\times  \acube{n}} \),  by
\begin{equation}
\label{eq:Psi_Cube}
\psi^{\oblong}_{m,n}\, \, = \sum_{\tau \in \shuf{m}{n}} \epsilon(\tau)\,  \tau*,
\end{equation}
where \(  \tau*\mbt := \left(t_{\tau^{-1}(1)}, \ldots, t_{\tau^{-1}(m+n)}\right) \). 
These maps   are related as follows. 
\begin{lemma}
\label{lem:Psi_mn}
For all \(  m, n \geq 0 \) the following diagram commutes. 
\begin{equation}
\label{eq:shuffle2}
\xymatrix{
\quad \asimplex{m+n} 
\ar@/_{1pc}/[d]^{\ \ \cong}_{\phi_{m+n}}
 \ar[rr]^{\psi^\Delta_{m,n}} & & \asimplex{m}\times \asimplex{n} \ar@/_{1pc}/[d]^{\ \ \cong}_{\phi_{m}\times \phi_n} \\
\quad  \acube{m+n} \ar@/_{1pc}/[u]_{\dif{m+n}} \ar[rr]_{\psi^\oblong_{m,n}} & & \acube{m}\times  \acube{n}
\ar@/_{1pc}/[u]_{\dif{m}\times \dif{n}}
}
\end{equation}
\end{lemma}

Let $\Gamma_\tau \subset \bbp^{m+n} \times \bbp^m \times \bbp^n$ be the closure of the graph of $\lambda_\tau$ \eqref{eq:lambda}. 
The projections 
\[
\begin{tikzcd}[row sep=small]
\bbp^{m+n} & &  \Gamma_\sigma \ar[ll, "\rho_{m,n}"'] \ar[rr, "\pi_{m,n}"] &&  \bbp^m \times \bbp^n.
\end{tikzcd}
\]
are surjective and equidimensional over $\Delta^{m+n}$ and $\Delta^m \times \Delta^n,$ respectively, allowing us to directly apply the current transform formalism in \cite{MR4498559} to get the next result.

\begin{proposition}
\label{prop:inv_iso}
For a fixed $\tau \in \mathsf{Sh}(m,n)$ the current transforms under $\Gamma_\tau$ and its transpose $\Gamma_\tau^t$ induce inverse isomorphisms
\(
\begin{tikzcd}
\nor{k}{\bbp^{m+n}|\infty}   \ar[rr, "\Gamma_{\tau \#}", bend left=6]  & & 
\nor{k}{\bbp^m\times \bbp^n|\infty}\ar[ll, "\Gamma^\#_{\tau}", bend left=6].
\end{tikzcd}
\)
\end{proposition}

From now on, let us  denote by \( \psi_{m,n} \) the simplicial EZ map \(  \psi_{m,n}^\Delta\), to simplify notation.
\begin{corollary}
\label{cor:CE-relation}
The  Eilenberg-Zilber morphisms $\psi_{m,n}$ induce homomorphisms
\[
\psi_{m,n \#} \ \colon \ \scrN^k(\bbp^{m+n}|\infty) \longrightarrow
\scrN^k(\bbp^m\times \bbp^n |\infty)
\]
satisfying the following relations. For all $m, n, r$ one has:
\begin{enumerate}[1)]
\item $\psi_{m,n\#}\circ \delta_\# \ = \ ( \delta \times \mathbb{1}_n)_\# \circ \psi_{m-1,n\#} \ + \  (-1)^m(\mathbb{1}_m \times \delta ) \circ \psi_{m,n-1 \#} $
\item $\left( \psi_{m,n\#} \times \mathbb{1}_r \right) \circ \psi_{m+n,r\#}  \ = \ (\mathbb{1}_m \times \psi_{m,n\#} )\circ \psi_{m,n+r\#}$.
\end{enumerate}
The latter identity simply states that the following diagram commutes.
\[
\begin{tikzcd}[row sep=0.2cm]
  &&  \nor{k}{\bbp^{m+n}\times \bbp^r | \infty} \ar[drr, "(\psi_{m,n} \times \mathbb{1}_r)_\#" ] & & \\
\nor{k}{\bbp^{m+n+r}  | \infty } \ar[urr, "\psi_{m+n, r \#}"] \ar[drr, "\psi_{m, n+r \#}"'] && &  & \nor{k}{\bbp^{m}\times\bbp^n \times \bbp^r | \infty }\\
  &&    \nor{k}{\bbp^{m}\times \bbp^{n+r} | \infty } \ar[urr, "(\mathbb{1}_m\times \psi_{n,r})_\#"']  &&  
\end{tikzcd}
\]
\end{corollary}
 
The homomorphisms in the corollary above can be assembled into maps of total product complexes, that we call the \emph{extended Eilenberg-Zilber maps}. 

\begin{definition}
\label{def:psi-Pi}
Define $\Psi^1 \colon \scrN^r_\Pi(1|\infty) \to\scrN^r_\Pi(2|\infty) $ and $\Psi^2 \colon \scrN^r_\Pi(2|\infty) \to\scrN^r_\Pi(3|\infty) $
as follows. 
\begin{enumerate}[a.]
\item
Given 
 \(\balpha  = \sum_{c\geq 0} \alpha^r_c\ t^c \  \in \ \scrN^r_\Pi(1|\infty), \) with \( \alpha^r_c \in \scrN^{r+c}(\bbp^c|\infty), 
\)
define 
\begin{equation}
\label{eq:Psi1}
\Psi^1(\balpha)  := \sum_{a_1,a_2\geq 0} \psi_{a_1,a_2 \#} \left( \alpha^r_{a_1+a_2} \right) \ x_1^{a_1}x_2^{a_2} \ \in \  \scrN^r_\Pi(2|\infty).
\end{equation}
\item
Given 
\( \bbeta  = \sum_{b_1, b_2\geq 0} \beta^r_{b_1, b_2}\  x_1^{b_1}x_2^{b_2} \   \in \ \scrN^r_\Pi(2|\infty)\) with \( \beta^r_{b_1, b_2} \in \scrN^{r+b_1+b_2}(\bbp^{b_1}\times \bbp^{b_2}|\infty),\)
define 
\begin{equation}
\label{eq:Psi2}
\Psi^2(\bbeta) =  \sum_{c_1,c_2, c_3\geq 0}   \Psi^2(\bbeta)_{c_1, c_2, c_3}\  y_1^{c_1}y_2^{c_2}y_3^{c_3}  \ \in \  \scrN^r_\Pi(3 |\infty),
\end{equation}
where \(  \Psi^2(\bbeta)_{c_1, c_2, c_3}  :=  (\psi_{c_1,c_2}\times \mathbb{1}_{c_3})_\#  ( \beta^r_{c_1+c_2, c_3}  )  - ( \mathbb{1}_{c_1} \times \psi_{c_2, c_3})_\#(\beta^r_{c_1,c_2+c_3}).\)
\end{enumerate}
\end{definition}

\begin{proposition}
\label{prop:Psi12}
The extended Eilenberg-Zilber maps  
\[
\Psi^1 \colon \scrN^*_\Pi(1|\infty) \to\scrN^*_\Pi(2|\infty) \quad\quad \text{and} \quad  \quad 
\Psi^2 \colon \scrN^*_\Pi(2|\infty) \to\scrN^*_\Pi(3|\infty) 
\] 
are morphisms of complexes.  Furthermore, $\Psi^2\circ \Psi^1 = 0.$
\end{proposition}
\begin{proof}
This is a straightforward calculation left to the reader. 
\end{proof}

\begin{remark}
\label{rem:extend_Psi}
One may wonder about the possibility  of defining maps $\Psi^{k}_\# \colon \scrN^*_\Pi(k|\infty) \to \scrN^*_\Pi(k+1|\infty) $ such that 
$\Psi^{k+1}_\# \circ \Psi^{k}_\# = 0$ for all $k\geq 0$, yielding a differential bigraded  algebra with potentially interesting applications.  This is not needed here. 
\end{remark}

\subsection{A character of the permutation Hopf algebra}
\label{subsec:char}
\newcommand{\Rdual}[1]{\check{\bbr}^{#1}}

In order to understand  the behavior of the fundamental triple \( \bbf = (\simplex{}, \Theta, \bbw) \)  under the Eilenberg-Zilber maps, we take a small detour that leads to a character of the permutation Hopf algebra \( \symg{}\mathsf{Sym}\). 

Consider the lattice \( \bbz^N\subset \bbr^N\) with canonical basis \( \{ \mbe_1, \ldots, \mbe_N\}\) and coordinates \( \mbt = (t_1,\ldots, t_N)\), and let
\(\Rdual{N}\)  be the dual of \( \bbr^N\), with dual basis \( \{ \check{\mbe}_1, \ldots, \check{\mbe}_N\}\) and coordinates \( \mby = (y_1, \ldots, y_N)\). 
The \emph{lattice-point transform} \( \Phi(P)(\mbx) \) of  a convex polytope \( P \subset \bbr^{N}\) is the formal Laurent series in the variables
%\( \bbz\llb x_1^\pm, \ldots, x_N^\pm \rrb \) 
\(  x_1, \ldots, x_N  \) defined by 
\begin{equation}
\label{eq:LPT}
\Phi(P)(x_1, \ldots, x_N) \ = \ \sum_{I\in P\cap \bbz^N}\ x_1^{i_1} \cdots x_N^{i_N}, \quad  \quad  I = (i_1, \ldots, i_N);
\end{equation}
see \cite{MR1175519}.

\begin{example}
\label{exmp:CN}
Consider the convex cone 
\begin{equation}
\label{eq:CN}
\calc_N := \{ (t_1, \ldots, t_N) \in \bbr^N \mid t_1\geq 0,\ t_1+t_2 \geq 0, \ldots, t_1 + \cdots + t_N \geq 0 \}.
\end{equation}
The lattice-point transform of \( \calc_N\) is given by
\begin{align}
\Phi & (\calc_N)(\mbx)  = \sum_{I\in \calc_N\cap \bbz^N} \ x_1^{i_1} \cdots x_N^{i_N}  = 
\sum_{I\in \calc_N\cap \bbz^N} \left(  \frac{x_1}{x_2}\right)^{i_1}\
\left(  \frac{x_2}{x_3}\right)^{i_1+i_2} \cdots \ \ \left(  \frac{x_{N-1}}{x_N}\right)^{i_1+\cdots + i_{N-1}} 
\left( x_N\right) ^{i_1+\cdots + i_N}     \\
& 
\label{eq:exmpCN} = 
\sum_{J \in \bbz_{\geq 0}^N}\  \left(  \frac{x_1}{x_2}\right)^{j_1}\
\left(  \frac{x_2}{x_3}\right)^{j_2} \cdots \ \ \left(  \frac{x_{N-1}}{x_N}\right)^{j_{N-1}} 
\left( x_N\right) ^{j_N}   =
\left( \frac{1}{1- \frac{x_1}{x_2}}\right) \left( \frac{1}{1- \frac{x_2}{x_3}}\right)\cdots \left( \frac{1}{1- \frac{x_{N-1}}{x_N}}\right) \left( \frac{1}{1-x_N}\right)  \notag \\
&  = \ 
\frac{x_2}{x_2-  x_1} \cdot \frac{x_3}{x_3- x_2} \cdots  \frac{x_N}{x_N -x_{N-1}}\cdot  \frac{1}{1-x_N}\   = \  
(x_1 \cdots x_N) \cdot \frac{1}{x_1(x_2-x_1)\cdots (x_N-x_{N-1})(1-x_N)}. \notag
\end{align}
\end{example}
%\smallskip

\begin{definition}
\label{def:alpha}
Let \( \bbz[\mbx] \)  be the  polynomial ring on a countably infinite set of variables \( \mbx = \{ x_1, x_2, \ldots \} \)  and let  \( \bbq(\mbx) \) be the corresponding field of rational functions. Let \( \bbz\llb \mbx, \mbx^{-1} \rrb\) be the ring of formal Laurent series on \( \mbx\).
 Given \( N \in \bbn\), define
\begin{equation}
\label{eq:alpha}
{\chi}_ N(\mbx) \ =  \
\chi_N(x_1, \ldots, x_N) : = \frac{1}{x_1(x_2-x_1) \cdots (x_N - x_{N-1}) (1-x_N)} \ \in \ \bbq(\mbx).
\end{equation}
\end{definition}

During our initial attempts to prove Theorem \ref{thm:pf-theta}, we encountered the intriguing  \emph{shuffle relations} in Proposition \ref{prop:aguiar} below. Thanks to the invaluable expertise of Marcelo Aguiar, who kindly outlined a more conceptual proof, we can present the results in a clearer and more elegant manner. The formulation of the next theorem is also due to Aguiar. 

\begin{theorem}
\label{thm:charHopf}
The collection  \( \{ \chi_N(\mbx) \mid N \in \bbn\} \) defines a character \( \bchi \) of the permutation Hopf algebra \( \perHA\) with values in the field \( \bbq(\mbx)\), 
defined by the assignment
\[
\bchi \ \colon \ \scrF_\sigma \ \longmapsto \chi_m\left(x_{\sigma(1)}, \ldots, x_{\sigma(m)}\right), \quad \sigma \in \symg{m},
\]
where \( \{  \scrF_\sigma \mid \sigma \in \symg{m}, m\geq 0 \} \) is the \( \scrF \)-basis for \( \perHA\).
\end{theorem}
\begin{proof}
This follows immediately from the following combinatorial result and from definitions.
\end{proof}

\begin{proposition}
\label{prop:aguiar}
The  collection of rational functions  \(  \{ \chi_N(\mbx), N\geq 0\} \) satisfies the shuffle relations
\[
\sum_{\tau \in \shuf{m}{n}} \ \chi_{m+n}\left( x_{\tau^{-1}(1)}, \ldots, x_{\tau^{-1}(m+n)} \right) \ = \ \chi_m(x_1, \ldots, x_m)\, \chi_n(x_{m+1}, \ldots, x_{m+n}),
\]
for all \( m, n \geq 0 \), where \( \shuf{m}{n}\) is the set of shuffles of type \( (m,n) \).
\end{proposition}

\begin{proof}
%[Proof of Proposition \ref{prop:aguiar}]
Let \( \btr{N} \subset \Rdual{N} \) be the topological \( N\)-simplex defined by 
\[ \btr{N} := \{ \mby = (y_1, \ldots, y_N) \mid 0 \leq y_1 \leq y_2 \leq \cdots \leq y_N \leq 1 \} 
= \mathsf{Cvx\, Hull}\la \mathbb{f}_1, \ldots, \mathbb{f}_N, {0} \ra ,
\]
where \( \mathbb{f}_j := \check{\mbe}_j  + \cdots + \check{\mbe}_N   \in \Rdual{N},\  j=1, \ldots, N.\) In  coordinates, \( \mathbb{f}_1 = (1, 1, \ldots, 1) \)  and   use \( \mathbb{1}_N := \mathbb{f}_1  \in \Rdual{N}\)  as   a distinguished vertex in \( \btr{N}. \)
Now, let \( T_{\bone_N}(\btr{N}) \) be the ``angle of  \( \btr{N}\) at \( \bone_N\)'' (see \cite[p. 26]{MR1234037}), consisting of all  vectors based  at \( \bone_N\) and pointing in the direction of points in \( \btr{N}\).  Equivalently, 
\begin{align*}
 T_{\bone_N}(\btr{N}) & =
  \mathsf{Cvx\, Hull}\left\{  \overrightarrow{\mathbb{f}_2 - \mathbb{f}_1},  \overrightarrow{\mathbb{f}_3 - \mathbb{f}_2},  \ldots, 
  \overrightarrow{0 - \mathbb{f}_1}
 \right\}  =
  \mathsf{Cvx\, Hull}\left\{  \overrightarrow{ - \check{\mbe}_1},  \overrightarrow{ - \check{\mbe}_1 - \check{\mbe}_2},  \ldots, 
  \overrightarrow{ - \check{\mbe}_1 - \check{\mbe}_2-\cdots  - \check{\mbe}_N}  \right\},
\end{align*}
where \( \overrightarrow{\mbv} := \{ t \mbv \mid t \in \bbr_{\geq 0} \} \) denotes the ray based at \( 0\) through  \( \mbv \neq 0 \in \Rdual{N} \).

The \emph{dual cone} \cite[p. 26]{MR1234037} to \( -  T_{\bone_N}(\btr{N}) \) is given by 
\begin{align}
\sigma_{\bone_N}\left(  - T_{\bone_N}(\btr{N}) \right) & := \ \left\{ \mbt \in \bbr^N \mid  - \lambda (\mbt) \geq 0, \ \text{ for all } \ \lambda \in T_{\bone_N}(\btr{N}) \right)\\
& =
\{ (t_1, \ldots, t_N) \in \bbr^N \mid t_1 \geq 0,\ t_1+t_2 \geq 0, \ldots,\ t_1 + t_2 + \cdots + t_N \geq 0 \} \notag  = \calc_N; \ \text{ see \ref{eq:CN}}. \notag
\end{align}

Now, consider the convex subset \( \btr{m}\times \btr{n} \subset \Rdual{m}\times \Rdual{n} \), which contains \( \bone_{m+n} = \bone_m\times \bone_n\) as a vertex. The  cone \( T_{\bone_{m+n}}(\btr{m}\times \btr{n})\) is the product cone \( T_{\bone_m} (\btr{m})\times T_{\bone_n}(\btr{n})\). As a result, 
the dual  to the opposite cone  \( - T_{\bone_{m+n}}(\btr{m}\times \btr{n})\) is given by
\begin{align}
\label{eq:prod_cone}
\sigma_{\bone_{m+n}}  \left(  - T_{\bone_{m+n}}(\btr{m}\times \btr{n}) \right)   & = 
\left\{
\mbt \in \bbr^{m+n} \left|  \ 
\begin{minipage}{9cm}
$t_1 \geq 0,\ t_1+t_2 \geq 0, \ldots,\ t_1 + t_2 + \cdots + t_m \geq 0, \ \text{ and }$ \\
$t_{m+1} \geq 0,\ t_{m+1}+t_{m+2} \geq 0, \ldots,\ t_{m+1} + t_{m+2} + \cdots + t_{m+n} \geq 0 $
\end{minipage} 
\right. \right\} \notag \\
&  = \ \calc_m \times \calc_n. 
\end{align}

The EZ construction in this context uses the action of  \(\symg{m+n} \) to define for each \( \tau \in \shuf{m}{n} \) the map
\begin{align*}
\lambda_{\tau^{-1}} \colon \btr{m+n} &\longrightarrow  \btr{m}\times \btr{n}\\
(y_1, \ldots, y_{m+n}) & \longmapsto (y_{\tau(1)}, \ldots, y_{\tau(m)})\times (y_{\tau(m+1)}, \ldots, y_{\tau(m+n)} ), 
\end{align*}
so that the collection of simplices \(\left\{ \btr{m+n}_\tau := \lambda_{\tau^{-1}}(\btr{m+n}) \mid \tau \in \shuf{m}{n} \right\} \) gives a triangulation of \( \btr{m}\times \btr{n} \).
Note that \( \lambda_{\tau^{-1}}(\bone_{m+n}) = \bone_{m+n} \), for all \( \tau \in \shuf{m}{n} \). Hence, \( \bone_{m+n} \) is a vertex in each simplex \( \btr{m+n}_\tau \). 

It follows that the collection of cones \( \left\{ T_{\bone_{m+n}}(\btr{m+n}_\tau) \mid \tau \in \shuf{m}{n} \right\} \)  forms a subdivision of the cone
\( T_{\bone_{m+n}}\left( \btr{m}\times \btr{n}\right) \)  in the sense of \cite[\S 1.2]{MR1175519}. Hence,
 \( \left\{ - T_{\bone_{m+n}}(\btr{m+n}_\tau) \mid \tau \in \shuf{m}{n} \right\} \)
 is  a corresponding subdivision of the opposite cone 
\( - T_{\bone_{m+n}}\left( \btr{m}\times \btr{n}\right) \).
For each \( \tau \in \shuf{m}{n} \),   denote the corresponding dual cones of this subdivision by
\begin{align*}
\calc_{m+n}^\tau & := \sigma_{\bone_{m+n}}\left( - T_{\bone_{m+n}}(\btr{m+n}_\tau) \right)    =
\left\{
\mbt  \mid  t_{\tau^{-1}(1)}\geq 0, \, t_{\tau^{-1}(1)} + t_{\tau^{-1}(2)} \geq 0, \, \ldots \, , 
t_{\tau^{-1}(1)}+ \cdots + t_{\tau^{-1}(m+n)} \geq 0
\right\}.
\end{align*}

Using this latter subdivision, one concludes from \cite[Corollaire, p.77]{MR1175519} that 
\begin{equation}
\label{eq:BrionCor} 
\Phi(\calc_m\times \calc_n) (\mbx) = \sum_{\tau\in \shuf{m}{n}}\ \Phi(\calc_{m+n}^{\tau})(\mbx).
\end{equation}
On the other hand, 
\begin{align*}
\Phi(\calc^\tau_{m+n})(\mbx) \ \ & := \sum_{J\in \calc_{m+n}^\tau \cap \bbz^{m+n}}\ x_1^{j_1} \cdots x_{m+n}^{j_{m+n}} 
\  \ =
\sum_{I \in \calc_{m+n} \cap \bbz^{m+n}}\ x_1^{i_{\tau(1)}} \cdots x_{m+n}^{i_{\tau(m+n)}} \\
& =
\sum_{I \in \calc_{m+n} \cap \bbz^{m+n}}\  x_{\tau^{-1}(1)}^{i_1}\  \cdots\  x_{\tau^{-1}(m+n)}^{i_{m+n}} \   = 
\Phi(\calc_{m+n})\left( x_{\tau^{-1}(1)}, \ldots,  x_{\tau^{-1}(m+n)} \right),
\end{align*}
and this identity, along with  \eqref{eq:BrionCor}, shows that 
\begin{align}
\Phi(\calc_m \times \calc_n)(x_1, \ldots, x_{m+n})   & = \Phi(\calc_m)(x_1, \ldots, x_m) \Phi(\calc_n)(x_{m+1}, \ldots, x_{m+n})  \notag  \\
& 
\label{eq:sum2}=
\sum_{\tau\in\shuf{m}{n}}\ \Phi(\calc_{m+n})\left( x_{\tau^{-1}(1)}, \ldots,  x_{\tau^{-1}(m+n)} \right)
\end{align}

Now, observe that  \( \chi_N(x_1, \ldots, x_N) = \frac{1}{x_1\cdots x_N} \cdot \Phi(\calc_N)(x_1, \ldots, x_N)\) for all \( N \in \bbn\), according to Example~\ref{exmp:CN}  and Definition \ref{def:alpha}.  Therefore,  the statement in the proposition follows by multiplying the terms of  second identity in \eqref{eq:sum2} by the \( \symg{N}\)-invariant element \( \frac{1}{x_1\cdots x_N} \ \in \ \bbq(\mbx) \).

\end{proof}

\begin{remark}
\label{rem:Hopf}
The Hopf algebra of   quasi-symmetric functions \( \mathscr{Q} \mathsf{Sym}\)  carries a notable character given by the multiple Zeta values
\( \textsc{mzv}\ \colon\ \mathscr{Q} \mathsf{Sym} \to \bbr\), and this follows from Kontsevich's representation of \( \zeta(n_1, \ldots, n_\ell) \)  as an integral over the \( n \)-simplex \( \blacktriangle^n\) with \( n = n_1+ \cdots + n_\ell\), along with   appropriate shuffle relations and the EZ   triangulations. Since \( \mathscr{Q} \mathsf{Sym}\) is as a quotient of the permutation Hopf algebra \( \symg{}\mathsf{Sym}\),  
  one may ask if there is  a subalgebra \( \cals \subset \bbq(\mbx) \)  along with a  \( \bbq\)-algebra homomorphism \(\Phi \colon \cals \to \bbr \), arising from periods of mixed Tate motives, so that 
\( \text{Image}(\boldsymbol{\chi}) \subset \cals\) and   the diagram 

\[
\begin{tikzcd}
\symg{}\mathsf{Sym} \ar[r, "\boldsymbol{\chi}"] \ar[d, "\pi"'] & \cals  \ar[d, "\Phi"] \ar[r, hook] &    \bbq(\mbx) \\
\scrQ \mathsf{Sym} \ar[r, "\textsc{mzv}"'] & \bbr & 
\end{tikzcd}
\]
commutes, where \( \pi \colon \symg{}\mathsf{Sym} \to\scrQ \mathsf{Sym} \) is the quotient map.
\end{remark}

\subsection{The fundamental triple under the extended Eilenberg-Zilber map}
\label{subsec:FT-EZ}

The  restriction of the components  $\lambda_\tau \) in the  EZ map \eqref{eq:Psi_Delta}
to the topological simplices \( \simplex{n} \subset \Delta^n\), along  with the orientation given by the sign of the permutation,   give a triangulation of the product $\simplex{m}\times \simplex{n}$. In particular,  the  map on integral currents
\(
\psi_{m,n\#} \ \colon \  \intcur{m+n}{\bbp^{m+n}|\infty} \to \intcur{m+n}{\bbp^m\times \bbp^n|\infty}
\)
induced by \( \psi_{m,n} \) satisfies
\begin{equation}
    \label{eq:pf-delta}
    \psi_{m,n\#}  ( \simplex{m+n}) = \simplex{m}\times \simplex{n},\ \ \text{ for all } m, n\geq 0.
\end{equation}%
 
\begin{proposition}
\label{prop:pf-delta}
Under the extended Eilenberg-Zilber map 
\(
\Psi^1 \ \colon \ \scrI^*_\Pi(1|\infty) \longrightarrow \scrI^*_\Pi(2|\infty)
\)
the simplex element \( \simplex{}= \sum_{b\geq 0} \simplex{b} \, x^b \) satisfies  \ \ 
\(
%\begin{equation}
%\label{eq:prod_simplex}
\Psi^1( \simplex{} ) = \simplex{} \boxtimes \simplex{} \ = \ \sum_{m, n\geq 0 }(\simplex{m}\times \simplex{n}) x_1^{m}x_2^{n}.
%\end{equation}
\)
\end{proposition}
\begin{proof}
By Definition \ref{def:exterior}, 
 $\simplex{}\boxtimes \simplex{}= \sum_{m, n }(\simplex{m}\times \simplex{n}) x_1^{m}x_2^{n}$, 
since  $\simplex{} \in  \scrI^0_\Pi(1|\infty)$ lies in degree $r=0$. The proposition now follows from \eqref{eq:pf-delta} and \eqref{eq:Psi1}.
\end{proof}
In a similar fashion, the polar element \( \Theta\) satisfies the same property.

\begin{theorem}
\label{thm:pf-theta}
For all $m, n \geq 0$  one has 
\(
\psi_{m,n\#}  ( \Theta_{m+n}) = \Theta_{m}\times \Theta_{n}.
\)
It follows that, under the extended Eilenberg-Zilber  map $\Psi^1 \ \colon \  \scrN^*_\Pi(1|\infty) \to \scrN^*_\Pi(2|\infty)$
the polar element $\Theta = \sum_{b\geq 0} \Theta_b \ x^b \in \scrN^0_\Pi(1|\infty)$ satisfies 
\[
\Psi^1(\Theta) = \Theta\boxtimes \Theta \ = \ \sum_{m, n \geq 0 }(\Theta_{m}\times \Theta_{n})\ x_1^{m}x_2^{n}.
\]
\end{theorem}
Proposition \ref{prop:aguiar} is essential for the proof of this theorem, along with the following observation.
\begin{remark}
\label{rem:proj_cur}
Let \(  f \colon X \to Y\) is a smooth map between oriented manifolds. 
If \(  R \) is a locally integral current on \(  X \) and \(  \theta  \) is a \(  \call^1_{\text{loc}} \)-form on \(  Y \) such that \(  R \wedge f^* \theta \) is a normal current on \(  X \), then 
\(  f_\sharp\left( R \wedge f^*\theta \right) \) is a normal current on \(  Y \) and
\(
f_\sharp \left(R \wedge f^*\theta \right) = ( f_\sharp R ) \wedge \theta. 
\)
In particular, if
 \(  \xymatrix{ X \ar@/^{0.2pc}/[r]^{g} & Y  \ar@/^{0.2pc}/[l]^{f} } \) are  inverse diffeomorphisms  and \(  R\wedge \omega \) is a normal current on \(  X \), with \(  R \) as above and \(  \omega \) a \(  \call^1_\text{loc} \)-form on X, then
\begin{equation}
\label{eq:adj}
f_\sharp \left(R \wedge \omega \right) = ( f_\sharp R ) \wedge g^*\omega. 
\end{equation}
\end{remark}
\begin{proof}[Proof of Theorem \ref{thm:pf-theta}]
In order to compare \( \psi_{m,n\#}(\Theta_{m+n})  \) and \( \Theta_m\times \Theta_n\),  we push them forward using  the isomorphism 
\( \phi_m \times \phi_n \ \colon\  \asimplex{m}\times \asimplex{n} \xrightarrow{\  \cong \  } \acube{m}\times \acube{n} \) \ and compare the resulting currents in \( \acube{m} \times \acube{n} = \acube{m+n}\). 

Recall   \ that \( \Theta_N := \frac{(-1)^{N(N+1)/2}}{\twopii^N}[\theta_N]\) \  \eqref{eq:theta}\ \  and denote the evident projections by 
\[ 
\begin{tikzcd}[column sep=0.6cm]
\asimplex{m}  & \ar[l,  "p"']  \asimplex{m}\times \asimplex{n} \ar[r, "q"] &  \asimplex{n} & 
  \text{ and }   & 
 \acube{m} & \ar[l,  "\hat p"']  \acube{m}\times \acube{n} \ar[r, "\hat  q"] &  \acube{n},
\end{tikzcd}  \]
Then, using  \eqref{eq:phi},  Remark \ref{rem:proj_cur}, and   the definition of exterior product of currents \cite[4.1.8]{MR0257325}, one gets 
\begin{align}
 (\phi_m\times \phi_n)_\# & ( [\theta_m]\times [\theta_n] )  = (-1)^{mn} (\phi_m\times \phi_n)_\# [p^*\theta_m \wedge q^* \theta_n]   = (-1)^{mn} [(\dif{m}\times\dif{n})^*(p^*\theta_m \wedge q^* \theta_n)] \notag \\
\label{eq:prod}
 & =
 (-1)^{mn} [(\dif{m}^*(p^*\theta_m)  \wedge \dif{n}^*( q^* \theta_n)]    = 
  (-1)^{mn} [(\hat p^* \dif{m}^*(\theta_m)  \wedge   \hat q^*\dif{n}^*( \theta_n)].  
 \end{align}

Set \( t_0 = 0\) and \( t_{N+1}=1\), so that the pull-back of  \( \Theta_N\) under \(  \dif{N}  \) is expressed as 
{\small 
\begin{align*}
 & \dif{N}^*(\Theta_N) =     \sum_{r=0}^N (-1)^r \dif{N}^*\left(  \frac{dz_0}{z_0} \wedge \cdots \wedge \widehat{ \frac{dz_r}{z_r} } \wedge \cdots \wedge \frac{dz_N}{z_N}\right)   =    \sum_{r=0}^N (-1)^r \dif{N}^*\left( \frac{z_r}{z_0\cdots z_N}\,   dz_0  \wedge \cdots \wedge \widehat{  dz_r  } \wedge \cdots \wedge  dz_N \right) \\
&   =  \left\{ \frac{1}{t_1(t_2-t_1)\cdots ( 1-t_{N}) }\right\} \sum_{r=0}^N (-1)^r  \{ t_{r+1}-t_{r}\}\,  \,   dt_1   \wedge \cdots \wedge \widehat{(dt_{r+1} - dt_r )} \wedge \cdots \wedge  ( -d t_N)\\
&   = \chi_N(\mbt) \sum_{r=0}^N (-1)^r  \{ t_{r+1}-t_{r}\}\,    dt_1   \wedge \cdots \wedge dt_r\wedge \widehat{(dt_{r+1} - dt_r )} \wedge \cdots   \wedge (dt_N-dt_{N-1})\wedge  ( -d t_N)\\
& =  
(-1)^N \chi_N(t_1, \ldots, t_N) \,   dt_1 \wedge \cdots \wedge dt_N,
\end{align*}
}
where  the last identity follows from a straightforward calculation. 

It follows that 
\begin{align}
\label{eq:dbl-push}
&  (\phi_m\times \phi_n)_\#   ( [\theta_m]\times [\theta_n] )  =  (-1)^{mn} [(\hat p^* \dif{m}^*(\theta_m)  \wedge   \hat q^*\dif{n}^*( \theta_n)]    \\
 & =
 (-1)^{mn} [(\hat p^* \left\{ (-1)^m \chi_m(t_1, \ldots, t_m) dt_1 \wedge \cdots \wedge dt_m \right\} \wedge   \hat q^*\left\{(-1)^n  \chi_n(t_1, \ldots, t_n) dt_1 \wedge \cdots \wedge dt_n\right\} ] \notag\\
\label{eq:push1}
 & =
 (-1)^{mn+m+n} \chi_m(t_1, \ldots, t_m) \chi_n(t_{m+1}, \ldots, t_{m+n})\ \ dt_1 \wedge \cdots \wedge dt_{m+n}
 \end{align}

On the other hand, from Lemma \ref{lem:Psi_mn} and Remark \ref{rem:proj_cur} one gets
\begin{align}
& (\phi_m\times \phi_n)_\# \left(  \psi_{m,n \#} [\theta_{m+n}]\right)   =  \psi_{m,n \#}^\oblong 
 \left( \phi_{m+n\#} [\theta_{m+n}] \right) = \psi_{m,n \#}^\oblong 
  [\dif{m+n}^* \theta_{m+n}]  \notag  \\
  & :=
  \sum_{\tau \in \shuf{m}{n}} \epsilon(\tau) \  \tau_\# \left(   [\dif{m+n}^* \theta_{m+n}]\right)   \ = \ 
\sum_{\tau \in \shuf{m}{n}} \epsilon(\tau) \  \tau_\# \left[   (-1)^{m+n} \chi_{m+n}(t_1, \ldots , t_{m+n})\ \  dt_1 \wedge \cdots \wedge dt_{m+n} \right] \notag  \\
& =
(-1)^{m+n} \left\{ \sum_{\tau \in \shuf{m}{n}} \chi_{m+n}\left(t_{\tau(1)}, \ldots, t_{\tau(m+n)}\right) \right\}\ \ dt_1 \wedge \cdots \wedge dt_{m+n} .\notag
\end{align}

The result now follows from the definition of \( \Theta_N\), along with \eqref{eq:dbl-push},  \eqref{eq:push1} and Proposition \ref{prop:aguiar}.
\end{proof}

\begin{corollary}
\label{cor:cycHCG}
The element $\bbw = \sum_{n\geq 0} W_n \ t^n \ \in \ \scrN^{-1}_\Pi(1|\infty)$ satisfies
\(
    D \Psi^1_\# ( \bbw ) =\Theta\boxtimes \Theta - \simplex{}\boxtimes \simplex{}.
\)
In particular, 
\[ \Psi^1 (\bbf) := ( \Psi^1 \simplex{},\ \Psi^1 \Theta,\ \Psi^1 \bbw) = (\simplex{}\boxtimes \simplex{},\ \Theta\boxtimes \Theta,\ \Psi^1 \bbw) \ \in \ \bbz^0_\Pi(2 | \infty ) 
\]
 is a cycle in the complex \( \bbz^*_\Pi(2 |\infty) \).
\end{corollary}
\begin{proof}
This follows   from the theorem, along with Propositions \ref{prop:Psi12} and \ref{prop:pf-delta}.
% and \ref{prop:univ_cycle_infty}.
\end{proof}

%%%%%%%%%%%%%.  \input{Prod-Reg.tex}
% !TEX root =  main_HCG_prod.tex
\section{Products and regulators}

In this section, we establish the desired multiplicative properties of our regulator map. The key result, Theorem \ref{thm:MAIN}, follows directly from the functorial properties exhibited by the transform of currents under equidimensional cycles, the behavior of the fundamental triple \( \bbf \)  under the extended EZ map, and the identification of the   classes of \(  \Psi^1(\bbf) \) and 
\( \bbf \hboxt \bbf\)  in \(  H^0 ( \bbz_\Pi^*(2 | \infty)) \).  This comparison is based on the calculation of the cohomology of the complex \( \scrN^*_\Pi(2\, | \infty) \), provided below.

\begin{theorem}
\label{thm:coh-2}
The  product complex $\scrN^*_\Pi(2|\infty)$ is quasi-isomorphic to $\underline{\bbc}$ (concentrated at $0$). In particular,
$H^{-1}(\scrN^*_\Pi(2|\infty)= 0.$
\end{theorem}
\begin{proof}
See Appendix \ref{subsec:coh}.
\end{proof}

\begin{corollary}
\label{cor:homology}
The elements 
\( 
\Psi^1(\bbf ) = (\simplex{}\boxtimes \simplex{}, \Theta \boxtimes \Theta, \Psi^1(\bbw) ) \)\ and \ \( \bbf \hboxt \bbf  = (\simplex{}\boxtimes \simplex{}, \Theta \boxtimes \Theta, \msfP_t(\bbf, \bbf) ) 
\)
represent the same homology class in \(  H^0 ( \bbz_\Pi^*(2 | \infty)) \). 
\end{corollary}
\begin{proof}
Consider the difference   \( \Psi^1(\bbf ) - \bbf \hboxt \bbf  =  \left( 0, 0, \Psi^1(\bbw) - \msfP_t(\bbf, \bbf)\right)  \). The boundary formula for \( ( \bbz_\Pi^*(2|\infty), \hat D) \) shows that 
\( \Psi^1(\bbw) - \msfP_t(\bbf, \bbf) \ \in \ \scrN^{-1}_\Pi(2|\infty) \) is a cycle in the complex \(  ( \scrN^{*}_\Pi(2|\infty), D). \)
Therefore, one can use Theorem \ref{thm:coh-2}  to find \( \msfT(\bbf) \in \scrN^{-2}_\Pi(2|\infty) \) such that 
\( D \msfT(\bbf) =  \Psi^1(\bbw) - \msfP_t(\bbf, \bbf) \). It follows that 
\( \hat D ( 0, 0, \msfT(\bbf) ) = \Psi^1(\bbf) \ - \ \bbf\hboxt \bbf.\)
\end{proof} 

\begin{remark}
\label{rem:explicit}
Write the element \( T(\bbf) \), from  the proof of the corollary above, in the form \( \msfT(\bbf) := \sum_{m,n\geq 0} \ T_{m,n}\ x_1^m x_2^n\), \ as in Definition \ref{def:tot}, with \(  T_{m,n} \in \scrN^{m+n-2}(\bbp^m\times \bbp^n | \infty) \).  By \eqref{eq:D}
\begin{align*}
D\msfT(\bbf) & = d \msfT(\bbf) + (-1)^{-2} \hat \delta\msfT(\bbf) \\
& =
\sum_{m,n\geq 0} \{ dT_{m,n} +  (-1)^{m+n} (\delta\times 1)_\#(T_{m-1,n}) +  (-1)^n (1\times \delta)_\#(T_{m,n-1})\}x_1^mx_2^n.
\end{align*}
In particular, for each \( m, n \geq 0 \) one has
\begin{align}
\label{eq:DRmn}
 dT_{m,n} \ & + \ (-1)^{m+n} (\delta\times 1)_\#(T_{m-1,n}) \ + \ (-1)^n (1\times \delta)_\#(T_{m,n-1}) 
  = \ \psi_{m,n\#}(W_{m+n}) \ - \ \msfP_t(\bbf_m, \bbf_n).
\end{align}
\end{remark}

\subsection{The multiplicative  property}

Let \( U \) and \( V \) be smooth quasiprojective varieties, with   
 NCD compactifications 
 \(
 \begin{tikzcd}[column sep=0.8cm]
 U \ar[r, hook] & X & D=X-U \ar[l, hook']
 \end{tikzcd}
 \text{ and } 
 \begin{tikzcd}[column sep=0.8cm]
V \ar[r, hook] & Y & E=Y-V. \ar[l, hook']
 \end{tikzcd}
 \)
  
 The regulator map from \cite{MR4498559} is a map of complexes
 \(
 \reg_U \ \colon \ \sbcxeq{p}{*}{U} \to \bbz(p)^{2p-*}_\scrD(U),
 \)
 where  
 \[
 \bbz(p)^{*}_\scrD(U)  \ := \ \mathsf{Cone}\left\{  \scrI^*_\text{loc}(p)(U)  \oplus F^p\norlog{*}{X}{D}  \xrightarrow{\epsilon - \imath}   \lnor{*}{U}  \right\}[-1] 
 \]
 calculates the Deligne-Beilinson cohomology of \( U \).  Here,  \( \scrI_\text{loc}^k(p) := \bbz(p) \cdot \scrI_\text{loc}^k\) is the sheaf of locally integral currents of degree \( k \) with coefficients in \( \bbz(p) := \twopii^p\bbz\).

The main technical tool in the definition of \(
 \reg_U \ \colon \ \sbcxeq{p}{*}{U} \to \bbz(p)^{2p-*}_\scrD(U),
 \)
 is the current transform operation associated to equidimensional cycles. 
It is shown in    \cite[Cor. 2.5]{MR4498559} that an equidimensional cycle \( \Upsilon \) of codimension \( p \) in \( U\times \Delta^n \), induces   a morphism of shifted complexes
\begin{equation}
\label{eq:curr_transf}
\Upsilon^\vee \colon \nor{*}{\bbp^n|\infty} \longrightarrow \norlog{*+2(p-n)}{X}{D}, \quad    T \longmapsto \Upsilon^\vee_T,
\end{equation}
which is  compatible with the Hodge filtration and preserves integral currents. 

Denote \( c_{p,n} := (-1)^{\frac{n(n-1)}{2}}\twopii^p \). Then, using  the transform operation, define
\begin{equation}
\label{eq:reg}
\reg_U(\Upsilon) \ := \ c_{p,n} \Upsilon^\vee_{\bbf_n} = \left( c_{p,n} \Upsilon^\vee_{\simplex{n}},\, c_{p,n} \Upsilon^\vee_{\Theta_n},\ c_{p,n}\Upsilon^\vee_{\bbw_n}\right) \ \in \ \bbz(p)^{*}_\scrD(U).
\end{equation}
 It is  shown in \cite{MR4498559}   that \( \cald\, \reg_U(\Upsilon) = \reg_U(\partial \Upsilon) \), where \( \cald \) here is the differential in Deligne-Beilinson's complex, and \( \partial  \) is the differential in the higher Chow  complex.

% \subsubsection{Exterior products in HCG and transforms}
 
 The exterior product in higher Chow groups 
 \begin{equation}
 \label{eq:ext-HCG} 
 \bar \times \colon \sbcxeq{p}{m}{U}\otimes \sbcxeq{q }{n}{V} \to \sbcxeq{p+q}{m+n}{U\times V}
 \end{equation}
is defined as \( A \bar\times B :=   \psi_{m,n}^*(A \times B) \) in \cite[\S 5]{Blo-HCG},  where 
  \(  \psi_{m,n}^* = \sum_{\tau \in \shuf{m}{n}} \epsilon(\tau) (1 \times \lambda_{\tau})^*  \)\  is the pull-back  under the EZ map \( \psi_{m,n}\), induced by the isomorphisms \( (1\times \lambda_\tau) \colon U\times V \times \Delta^{m+n} \xrightarrow{\cong} (U\times \Delta^m) \times (V \times \Delta^n) \), and
   \( A \times B \) is the product cycle in \( (U\times \Delta^m) \times (V \times \Delta^n) \).
   
 On the other hand, to obtain the exterior product in Deligne-Beilinson cohomology one  defines, for    \( t \in \bbc\),  a pairing of complexes  suitably adapted from  \cite{MR944991} or \cite[Defn. 4]{MR2818720}  to complexes of currents 
 \begin{equation}
 \label{eq:extDBcoh}
 \boxtimes_t \ \colon \ \bbz(p)^*_\scrD(U) \otimes \bbz(q)^*_\scrD(V) \longrightarrow \bbz(p+q)^*_\scrD(U\times V),
 \end{equation}
 in a fashion similar to the constructions in  \S\ref{subsec:pairings}. 
 First, given    \( R = (R_1, R_2, R_3) \in \ \bbz(p)^m_\scrD(U) \) and
 \( S = (S_1, S_2, S_3) \in \bbz(q)^n_\scrD(V) \), define
 \begin{equation}
 \label{eq:del_prod}
 R\, \hboxt S := (R_1 \times S_1, \, R_2 \times S_2, \, \msfP_t(R,S) ),
 \end{equation}
 where \( \msfP_t(R,S) ) = R_3 \times \Xi_t(S_1,S_2)  +   \Xi_{(1-t)}(R_1,R_2)\times S_3 \). As in  Proposition \ref{prop:cupt}, one shows that \( \boxtimes_t\) is a map of complexes, whose homotopy class does not depend on \( t \in \bbc\). 
Now, for \( R\) and \( S\) as above, define
\begin{equation}
\label{eq:boxt}
R \boxtimes_t  S \ := \ (-1)^{mn} R\ \hboxt \ S.
\end{equation} 
\begin{remark}
\label{rem:cup}
With this sign convention, the \emph{cup product}  
\[
\cup \ \colon \ H^k_\scrD(U;\bbz(p)) \otimes H^\ell_\scrD(U; \bbz(p)) \longrightarrow H^{k+\ell}_\scrD(U;\bbz(p+q))
\]
is defined on classes represented by \( R\) and \( S \)  as
\( [R]\cup [S] \ = \ \gamma^*\left( [R \boxtimes_t S] \right) \), where \( \gamma \colon U \to U \times U \) is the diagonal map.  
This is   compatible with the forgetful functor into singular cohomology; see \cite[p. 328]{MR1700700}.
\end{remark}
 
The behavior of current transforms under products of algebraic cycles is explained below.

 \begin{lemma}
 \label{lem:taut}
 Given  currents \( R \in \nor{k}{\bbp^m|\infty},\ S \in \nor{\ell}{\bbp^n|\infty} \) and \( T \in \nor{d}{\bbp^{m+n}|\infty} \), and equidimensional cycles  \( A \in \sbcxeq{p}{m}{U} \) and \( B \in \sbcxeq{q}{n}{V} \), the following holds.
 \begin{enumerate}[i.]
 \item \label{it:taut1} \( (A \times B )^\vee_{R\times S} \ = \ A^\vee_T \times B^\vee_S\).
 \item \label{it:taut2} \( (A\bar \times B)^\vee_T \ = \ (A\times B)^\vee_{\psi_{m,n\#}(T)} \).
 \end{enumerate}
 \end{lemma}
 \begin{proof}
The result follows from properties of  current transforms shown in \cite{MR4498559} and from the definitions. For the reader's convenience,   details are written in Appendix \ref{subsec:PfLemma1}.
 \end{proof}

 \begin{theorem}
\label{thm:MAIN}
Let \( U \) and \( V\) be smooth quasi-projective varieties over \( \bbc\). Then, for any fixed \( t \in \bbc\), the following diagram of morphism of complexes commutes up to natural homotopy
\begin{equation}
\label{eq:up-to}
\begin{tikzcd}
\sbcxeq{p}{*}{U} \otimes \sbcxeq{q}{*}{V} \ar[rr,"\bar \times"] \ar[d, "\reg_U\otimes\,  \reg_V"']  & &  \sbcxeq{p+q}{*}{U\times V} 
\ar[d, "\reg_{U\times V}"] \\
\dcx{p}{2p-*}{U}\otimes \dcx{q}{2q-*}{V} \ar[rr,"\boxtimes_t"'] & &  \dcx{p+q}{2(p+q)-*}{U\times V},
\end{tikzcd}
\end{equation}
where 
\( \bar \times \) and \( \boxtimes_t \) are the exterior products in the equidimensional higher Chow complexes and Deligne-Beilinson current complexes, respectively.
\end{theorem}
\begin{proof}
In the definition of the regulator map \eqref{eq:reg},  the current transforms are normalized with the constants \( c_{d,N} := (-1)^{\frac{N(N-1)}{2}}\twopii^d\), which satisfy 
\begin{equation}
\label{eq:cpq}
c_{d,N} = (-1)^{N-1}c_{d,N-1} \quad \text{ and } \quad c_{p+q,m+n} =   (-1)^{mn}c_{p,m}c_{q,n}.  
 \end{equation}
Now,    using the elements   \( T_{m,n} \in \scrN^{m+n-2}(\bbp^m\times \bbp^n | \infty) \)   from  Remark~\ref{rem:explicit},    define homomorphisms
\begin{equation*}
H^d \colon  \left\{ \sbcxeq{p}{*}{U} \otimes \sbcxeq{q}{*}{V} \right\}^d := \bigoplus_{m+n = d} \ \sbcxeq{p}{m}{U} \otimes \sbcxeq{q}{n}{V}  \longrightarrow \bbz(p+q)^{2(p+q)+1-d}_\scrD(U\times V) 
\end{equation*}
by sending
%\begin{equation}
%\label{eq:AB}
 \(  A \in \sbcxeq{p}{m}{U}\)  and  \(B \in  \sbcxeq{q}{n}{V} \)
% \end{equation}
to \( H^{d} (A\otimes B) :=  \left( 0, 0, c_{p+q, m+n} \ (A\times B)^\vee_{T_{m,n}}\right) .\)
We will show  that  \( H^*\) defines a homotopy between 
\( \left( \reg_{U\times V}\right) \circ \bar \times \ \)  and \(\ \boxtimes_t \circ \left( \reg_U \otimes \reg_V\right).\) 

\begin{claim}
\label{claim:tech}
The homomorphisms \( H^d \) are natural on  pairs \( (U, V) \). Furthermore, 
if \( A, B \) are as above, then
\[
\cald \circ H^{m+n}(A\bar \times B) \ - \ H^{m+n-1}\circ \partial (A\bar \times B)  \ = \  \reg_{U\times V}(A\bar \times B)  \ - \ \reg_U(A) \boxtimes_t \reg_V(B),
\]
where \( \cald \) is the differential in the Deligne-Beilinson complex \( \bbz(p+q)^*_\scrD(U\times V) \), and \( \partial \) is the differential in the product complex \( \sbcxeq{p}{*}{U} \otimes \sbcxeq{q}{*}{V}\).
\end{claim}
The naturality assertion in the Claim follows from the  properties of the current transforms proven in \cite{MR4498559}. Now,
using Lemma \ref{lem:taut} one obtains
\begin{equation}
\begin{alignedat}{4}
& (A\bar \times B)^\vee_{\bbf_{m+n}} && : =
  \left( (A\bar \times B)^\vee_{\simplex{m+n}},
(A\bar \times B)^\vee_{\Theta_{m+n}}, (A\bar \times B)^\vee_{W_{m+n}} \right)    && \\
&  &&  = 
  \left( (A  \times B)^\vee_{\psi_{m,n\#}(\simplex{m+n})},
(A \times B)^\vee_{\psi_{m,n\#}(\Theta_{m+n})}, (A  \times B)^\vee_{\psi_{m,n\#}(W_{m+n})} \right)     && \text{ \quad (by Lem. \ref{lem:taut}.\ref{it:taut2})} \\
\label{eq:bar_prod}
&  && =  
  \left( (A  \times B)^\vee_{\simplex{m}\times \simplex{n}},\ 
(A \times B)^\vee_{ \Theta_{m }\times \Theta_n}, \ (A  \times B)^\vee_{\psi_{m,n\#}(W_{m+n})} \right)  && \text{ \quad (see \S\,\ref{subsec:FT-EZ})}  \\
& && 
=   \left( A^\vee_{\simplex{m}}  \times B^\vee_{  \simplex{n}},\ 
A^\vee_{ \Theta_m} \times B^\vee_{  \Theta_n},\  (A  \times B)^\vee_{\psi_{m,n\#}(W_{m+n})} \right)    && \text{ \quad (by Lem. \ref{lem:taut}.\ref{it:taut1}) } \\
&  &&
=
\left( A^\vee_{\simplex{m}}  \times B^\vee_{  \simplex{n}},
A^\vee_{ \Theta_m} \times B^\vee_{  \Theta_n}, \msfP_t(A^\vee_{\bbf_m}, B^\vee_{\bbf_n}) \right)   \\
&  && + \left( 0, \ 0,\  (A  \times B)^\vee_{\psi_{m,n\#}(W_{m+n})} - \msfP_t(A^\vee_{\bbf_m}, B^\vee_{\bbf_n}) \right)   && \\
& && =  A^\vee_{\bbf_m} \hboxt B^\vee_{\bbf_n} \ + \  \left( 0, 0, (A  \times B)^\vee_{\psi_{m,n\#}(W_{m+n})} - \msfP_t(A^\vee_{\bbf_m}, B^\vee_{\bbf_n})\right).   && \text{\quad (by \eqref{eq:del_prod})}
\end{alignedat}
\end{equation}

Since the operation \( \Xi_t(R,S) \) is a linear combination of the currents \( R, S \), one has
\begin{equation}
\begin{alignedat}{4}
& \msfP_t(A^\vee_{\bbf_m}, B^\vee_{\bbf_n}) && := 
A^\vee_{W_m} \times \ \Xi_t(B^\vee_{\simplex{n}},B^\vee_{\Theta_n}) \ + \  \Xi_{(1-t)}(A^\vee_{\simplex{n}},A^\vee_{\Theta_n})\times B^\vee_{W_n} &&  \\
&  && =
A^\vee_{W_m} \times \ B^\vee_{\Xi_t(\simplex{n},\Theta_n)} \ + \  A^\vee_{\Xi_{(1-t)}(\simplex{n},\Theta_n)} \times B^\vee_{W_n} && \\
\label{eq:fromPt}
& &&   = 
 (A \times B)^\vee_{W_m \times  \Xi_t(\simplex{n}, \Theta_n)} \ + \ (A \times B)^\vee_{\Xi_{(1-t)}(\simplex{n},\Theta_n) \times W_n} &&  \text{\quad (by Lem. \ref{lem:taut}.\ref{it:taut1})} \\
& && = 
(A\times B )^\vee _{W_m \times  \Xi_t(\simplex{n}, \Theta_n) \ + \ \Xi_{(1-t)}(\simplex{n},\Theta_n) \times W_n} \ =: \ (A \times B)^\vee _{\msfP_t(\bbf_m,\bbf_n)}. &&
\end{alignedat}
\end{equation}

Putting identities \eqref{eq:bar_prod}, \eqref{eq:fromPt} and \eqref{eq:Pt} together gives
\begin{align*}
(A\bar \times B)^\vee_{\bbf_{m+n}} & : =
A^\vee_{\bbf_m} \hboxt B^\vee_{\bbf_n} \ + \ \left(0, 0, (A  \times B)^\vee_{\psi_{m,n\#}(W_{m+n})} - (A \times B)^\vee _{\msfP_t(\bbf_m,\bbf_n)} \right) \\
& = 
A^\vee_{\bbf_m} \hboxt B^\vee_{\bbf_n} \ + \ \left(0, 0, (A  \times B)^\vee_{ \psi_{m,n\#}(W_{m+n})- \msfP_t(\bbf_m,\bbf_n) }\right),
\end{align*}
and it follows from \eqref{eq:DRmn} that 
\begin{align*}
& (A\bar \times B)^\vee_{\bbf_{m+n}}\  - \  
A^\vee_{\bbf_m} \hboxt B^\vee_{\bbf_n}    = \left(0, 0,\  (A  \times B)^\vee_{ dT_{m,n} +  (-1)^{m+n} (\delta\times 1)_\#(T_{m-1,n})  +  (-1)^n (1\times \delta)_\#(T_{m,n-1}) }\right) \\
& = 
\left(0, 0,\   d (A  \times B)^\vee_{T_{m,n}} \  +\   (-1)^{m+n} (A  \times B)^\vee_{(\delta\times 1)_\#(T_{m-1,n})}\ +\  (-1)^n (A  \times B)^\vee_{(1\times \delta)_\#(T_{m,n-1}) }\right) \\
& = 
\left(0, 0, \  d (A  \times B)^\vee_{T_{m,n}}  \ +\   (-1)^{m+n} (\partial A  \times B)^\vee_{T_{m-1,n}}  +  (-1)^n (A  \times \partial B)^\vee_{T_{m,n-1} }\right).
\end{align*}
Now, multiplying both sides of the equation above by \( c_{p+q,m+n} \) and using \eqref{eq:cpq} along with \eqref{eq:boxt} gives
\begin{align*}
& \reg_{U\times V}(A\bar \times B) \  - \ \reg_U(A)\boxtimes_t \reg_V(B) \\
& =
\left(0, 0,  c_{p+q,m+n}\ d (A  \times B)^\vee_{T_{m,n}}\right) 
\ +\  (-1)^{m+n}c_{p+q,m+n}\left( 0, 0,  \left\{ (\partial A  \times B)^\vee_{T_{m-1,n}}  +  (-1)^m (A  \times \partial B)^\vee_{T_{m,n-1} } \right\}\right) \\
 & =
 \cald  H^{m+n}(A  \bar \times B)  - c_{p+q,m+n-1} \left( 0, 0, \left\{ (\partial A  \times B)^\vee_{T_{m-1,n}}  +  (-1)^m (A  \times \partial B)^\vee_{T_{m,n-1} } \right\} \right) \\
 & =
  \cald  \circ H^{m+n}(A  \bar \times B)  \ - \ H^{m+n-1}\circ \partial (A\bar \times B)  .
\end{align*}

This concludes the proof of Claim \ref{claim:tech}, and the theorem follows. 

\end{proof} 

\begin{corollary}
\label{cor:ext}
If \( U \) and \( V\) are smooth quasi-projective varieties over \( \bbc\), then the regulator maps 
%\[ \reg_U \colon \shcg{p}{r}{U}\to \dcz{2p-r}{p}{U}\] and \( \reg_V \colon \shcg{q}{s}{V}\to \dcz{2q-s}{q}{V}\) 
commute with exterior product. In other words,
\[ \reg_{U\times V}(\alpha \times_{\text{\tiny CH}} \beta)  = \reg_U(\alpha) \times_\scrD \reg_V(\beta)\  \in\  \dcz{2(p+q)-(r+s)}{p+q}{U\times V},\] where 
\(  \times_{\text{\tiny CH}} \) and \( \times_\scrD\) are the exterior products in higher Chow groups and Deligne-Beilinson cohomology, respectively.
\end{corollary}

For smooth varieties, one has a natural isomorphism \( \shcg{p}{n}{U}\cong H_\scrM^{2p-n}(U;\bbz(p)) \) between higher Chow groups and Voevodsky's motivic cohomology  \cite[Lect. 19]{MR2242284}, which gives our last result.

\begin{corollary}
\label{cor:cup}
The regulator map \( \reg_U \colon  H_\scrM^{*}(U;\bbz(\bullet)) \to  \dcz{*}{\bullet}{U}\) is a homomorphism of bigraded rings, and is a natural transformation of functors in the category of smooth complex quasiprojective varieties.
\end{corollary}
\begin{proof}
This follows from definitions and Remark \ref{rem:cup}.
\end{proof}

%%%%%%%%%%%%%.  \input{Appendix}
% !TEX root =  main_HCG_prod.tex

\appendix

\section{Technical proofs}

\subsection{Proof of Theorem \ref{thm:coh-2}}
\label{subsec:coh}
% \lipsum[1]
Given $m, n\geq 0,$ consider the correspondence 
\(
C_{m,n}\subset \bbp^m \times \bbp^n \times \bbp^1 \times \bbp^m \times \bbp^n
\)
defined as the locus of those points
$([\mbu], [\mbv], [s_0:s_1], [\mbx], [\mby])$ satisfying
\begin{align*}
    x_i \ve(\mbu) & = \{ s_1 u_i + s_0\} \ve(\mbx), \ i= 0, \ldots, m \quad \text{ and } \quad  y_j \ve (\mbv)  = \{ s_1 v_j + s_0\} \ve(\mby), \ j=0, \ldots, n,
\end{align*}
where $\ve(\mbu) = u_0 + \cdots + \ u_m $\  and\ \( \ve(\mbv) = v_0+\cdots +v_n\).
It is easy to see that   $C_{m,n}$ is the closure of the graph of the linear homotopy that contracts $\Delta^m \times \Delta^n$ into the barycenter $(B_m, B_n)$. Hence, it is equidimensional and dominant over \( \bbp^m \times \bbp^n \times \bbp^1 \).

 It follows from \cite[\S 2]{MR4498559} that $C_{m,n}$ induces a transform homomorphism
 \begin{equation}
     C^\vee_{m,n} \colon \nor{k}{\bbp^m\times \bbp^n\times \bbp^1|\infty} \longrightarrow \nor{k-2}{\bbp^m\times \bbp^n |\infty}.
 \end{equation}

\begin{definition}
\label{def:Cmn}
For $k\geq 1$, define $C_{m,n\sharp}\ \colon\ \nor{k}{\bbp^m\times \bbp^n |\infty} \to \nor{k-1}{\bbp^m\times \bbp^n |\infty}$
as the composition
\begin{equation}
\label{eq:compo}
    \begin{tikzcd}
   \nor{k}{\bbp^m\times \bbp^n |\infty} 
   \ar[r, "{-}\times {\simplex{\, 1}}"]  \ar[rr, "C_{m,n\sharp}"', bend right=10]    &  \nor{k+1}{\bbp^m\times \bbp^n \times \bbp^1 |\infty} \ar[r, "C^\vee_{m,n}"]
   &  \nor{k-1}{\bbp^m\times \bbp^n |\infty}.
    \end{tikzcd}
\end{equation}
\end{definition}

\begin{lemma}
\label{lem:htpy}
Given $T\in \nor{k}{\bbp^m\times \bbp^n|\infty}$, with $0 \leq k \leq 2(m+n)-1$, one has
$T = dC_{m,n\sharp}(T) + C_{m,n\sharp}(dT).$
\end{lemma}
\begin{proof}
It follows from the definitions and the properties shown in \cite[\S 2]{MR4498559} that
\begin{align*}
    dC_{m,n\sharp}(T) & = dC^\vee_{m,n}(T \times \simplex{1}) = 
    C^\vee_{m,n}\{d(T \times \simplex{1})\}   =
    C^\vee_{m,n}\left(-dT \times \simplex{1} +  [T\times \mbe_1 - T\times \mbe_0] \right) \\
    & =
    - C^\vee_{m,n}\left(dT \times \simplex{1}\right)  +  C^\vee_{m,n}(T\times \mbe_1)  - C^\vee_{m,n}(T\times \mbe_0)   =
    - C_{m,n\sharp }\left(dT \right) + T - C^\vee_{m,n}(T\times \mbe_0).
\end{align*}
From the definition of $C_{m,n}$ and the fact that $T$ vanishes suitably at infinity (see  \S \ref{subsec:interest}), one sees that the support of $C^\vee_{m,n}(T\times \mbe_0)$ is contained in the singleton $\{ (B_m, B_n) \}$. From the \emph{support theorem} for normal currents \cite[\S 4.1.15]{MR0257325}, it folows that if $\dim(T)>0$ (i.e. $\deg(T) < 2m+2n$) then $C^\vee_{m,n}(T \times \mbe_0) = 0. $ Therefore, for $T \in \nor{k}{\bbp^m\times \bbp^n|\infty}$ with $0 \leq k < 2m+2n$ one has
$d C_{m,n\sharp}(T) + C_{m,n\sharp}(dT) = T.$
\end{proof}
\begin{corollary}
Given $\mbb \in \bbz^2_{\geq 0}$, the cohomology groups of the complex $\nor{*}{\bbp(\mbb)|\infty}$ computes the cohomology with compact supports of $\Delta^{b_1}\times \Delta^{b_2},$ that is, 
\[
H^r(\nor{*}{\bbp(\mbb)|\infty} \cong
H_c^r( \Delta^{b_1}\times \Delta^{b_2}; \bbc) \cong\begin{cases}
\bbc, & r = 2|\mbb| \\
0, & \text{otherwise}.
\end{cases}
\]
\end{corollary}

\begin{proof}[Proof of Theorem \ref{thm:coh-2}]
Given $k\geq 1$ and $p\geq 0$, define a filtration on $\scrN^*_\Pi(k|\infty)$ by setting
\(
F^p\scrN^r_\Pi(k|\infty) = \prod_{|\mbb|\geq p} \nor{r+|\mbb|}{\bbp(\mbb | \infty)}
\)
to be the subgroup of $\scrN^r_\Pi(k|\infty)$ consisting of the elements whose coordinates in 
$\scrN^{r+|\mbb|}(\bbp(\mbb)|\infty)$ are zero if $|\mbb|<p.$
Therefore, one obtains a filtered complex
\[
\scrN^*_\Pi(k|\infty)= F^0\scrN^*_\Pi(k|\infty) \supset \cdots \supset F^p\scrN^*_\Pi(k|\infty) \supset \cdots 
\]
such that
\(
%\begin{equation}
%    \label{eq:Q1}
\grF{p}\left\{ \scrN^*_\Pi(k|\infty) \right\} :=    F^p\scrN^*_\Pi(k|\infty)/F^{p+1}\scrN^*_\Pi(k|\infty) \ \cong \ \prod_{|\mbb| = p} \nor{r+p}{\bbp(\mbb)|\infty}.
%\end{equation}
\)
The result is a spectral sequence $E_\bullet^{*,*}(k)$ with
\begin{align}
E_1^{p,q}(k) & = H^{p+q}\left(  F^p\scrN^*_\Pi(k|\infty)/F^{p+1}\scrN^*_\Pi(k|\infty) \right) = \prod_{|\mbb|=p} H^{p+q} \left( \nor{*+p}{\bbp(\mbb)|\infty} \right) \\
& \cong
\prod_{|\mbb|=p} H_c^{2p+q} \left( \bbp(\mbb) - H_\infty(\mbb); \bbc\right) = \prod_{|\mbb|=p} H_c^{2p+q} \left( \Delta^{b_1} \times \cdots \times \Delta^{b_r}; \bbc\right) \ \cong 
\begin{cases}
\prod_{|\mbb| = p} \bbc &, q=0 \\
0 &, \text{otherwise}.
\end{cases} \notag
\end{align}
It is convenient to express the calculation above as
\begin{equation}
    \label{eq:sseq}
    E_1^{p,q}(k) \cong \begin{cases}
    \bbc[x_1, \ldots, x_k]_p &, q=0\\
    0 &, \text{otherwise},
    \end{cases}
\end{equation}
where $\bbc[x_1, \ldots, x_k]_p$ denotes the homogeneous polynomials of degree $p$ in the variables $x_1, \ldots, x_k.$ In this description, $x_1^{b_1} \cdots x_k^{b_k}$ corresponds to the generator of
$H^{2p}_c(\Delta^{b_1} \times \cdots \times \Delta^{b_r}; \bbc) $, which can be given by the $0$ dimensional current $\eta_{\mbb} \in \scrN^{2p}(\bbp(\mbb)|\infty)$ represented by integration over  $B(\mbb) = (B_{b_1}, \ldots, B_{b_k})$, where $B_b = [1:1:\cdots:1] \in \bbp^b$ is the barycenter.

By definition, $\hat \delta(x_1^{b_1} \cdots x_k^{b_k}) = \sum_{j=1}^k (-1)^{b_j + \cdots + b_k} (1\times \cdots \times \delta \times \cdots \times 1)_\#(x_1^{b_1} \cdots x_k^{b_k})$; see \eqref{eq:dalpha}.
On the other hand, it is clear that the differential
$\delta \colon \scrN_0(\bbp^p|\infty) \to \scrN_0(\bbp^{b+1}|\infty)$ is $0$ when $b$ is even and an isomorphism when $b $ is odd. Therefore, 
$$\delta^j(x_1^{b_1} \cdots x_k^{b_k})=
\begin{cases}
0 &, \text{ if } b_j \text{ is even},\\
x_1^{b_1} \cdots b_j^{b_j+1}\cdots x_k^{b_k} &, \text{ if } b_j \text{ is odd}.
\end{cases}$$
This suffices to calculate the differentials in $E_1^{*,*}(k)$ for all $k$, but   the case $k=2$ suffices.

Write $E_1^{p,0}(2) = \bbc[x,y]_p$ and denote
$\tau(m) = \begin{cases}
0 &,  \text{ if $m$ is even,} \\
1 &, \text{ if $m$ is odd}.
\end{cases}
$

\noindent Then
%\begin{align*}
 \(   \hat\delta (x^ay^b)  = (-1)^{a+b}\delta_1(x^ay^b) + (-1)^b \delta_2(x^ay^b)   = 
(-1)^{a+b}\tau(a)\ x^{a+1}y^b + (-1)^b \tau(b) x^ay^{b+1}.\)
% \end{align*}
From now on, it becomes  easier to identify $\bbc[x,y]_p$ with $\bbc[x]_{\leq p}$ using $x\mapsto x$ and $y\mapsto 1$, and hence
$\hat\delta (x^a) = (-1)^p\tau(a)\ x^{a+1} + (-1)^{p-a}\tau(p-a)\ x^a$. \ 
Therefore, for $\balpha = \sum_{a=0}^p \alpha_a x^a$, one has
\begin{align*}
    \hat\delta \balpha & = 
    (-1)^p \sum_{a=0}^p \tau(a) \alpha_a x^{a+1} + (-1)^p \sum_{a=0}^p (-1)^a \tau_{p-a}\alpha_a x^a \\
    &  = 
    (-1)^p \tau(p) \alpha_p \cdot x^{p+1} + (-1)^p\tau(p) \alpha_0 \cdot 1   + (-1)^p \sum_{r=1}^p \left\{ \tau(r-1)\ \alpha_{r-1} + (-1)^r  \tau_{p-r}\ \alpha_r \right\} x^r .
\end{align*}

If  $p=2r$ is even, it follows that
$ \hat \delta \left(\sum_{a=0}^{2r} \alpha_a x^a\right) = \alpha_1 (x^2-x^1) + \alpha_3(x^4 - x^3) + \cdots + \alpha_{2r-1}(x^{2r}- x^{2r-1}).$ Therefore
\begin{equation}
    \label{eq:p_even}
    \begin{cases}
    \ker(\hat\delta^{2r}) & = \bbc \langle 1, x^2, \ldots , x^{2r} \ra = \text{ polynomials of even degree}\\
    \operatorname{Im}(\hat \delta^{2r}) & = \bbc \langle x^2-x, \ldots , x^{2r}-x^{2r-1} \ra .
    \end{cases}
\end{equation}

Similarly, when $p=2r+1$ is odd one obtains
\begin{align*}
\hat \delta \balpha & = \hat \delta \left(\sum_{a=0}^{2r+1} \alpha_a x^a\right)  = - \alpha_{2r+1} x^{2r+2} - \alpha_0 \cdot 1 - \left\{ (\alpha_1 + \alpha_2)x^2 + (\alpha_3 + \alpha_4)x^4 + \cdots + (\alpha_{2r-1}+\alpha_{2r})x^{2r}\right\}.
\end{align*}
Hence
\begin{equation}
    \label{eq:p_even2}
    \begin{cases}
    \ker(\hat\delta^{2r+1}) & = \bbc \langle x^2-x, \ldots , x^{2r}-x^{2r-1} \ra  \\
    \operatorname{Im}(\hat \delta^{2r+1}) & = \bbc \langle 1, x^2, \ldots , x^{2r} \ra = \text{ polynomials of even degree}.
    \end{cases}
\end{equation}
It is clear from these calculations that for $p\geq 1$ one obtains $\operatorname{Im}(\hat \delta^{p-1}) = \ker(\hat \delta^p).$ It now follows that
\begin{equation}
    \label{eq:E_infty}
    E_2^{p,q} \cong E_\infty^{p,q} \cong \begin{cases}
    \bbc, & p=q=0 \\ 
    0, & \text{otherwise}.
    \end{cases}
\end{equation}
This concludes the proof of Theorem \ref{thm:coh-2}.
\end{proof}

\subsection{Proof of Lemma \ref{lem:taut}}
\label{subsec:PfLemma1}
 \begin{proof}
Consider the  commutative diagram. 
 \[
 \begin{tikzcd}
 & &  U\times V  \ar[dll, bend right=10, "pr_1"'] \ar[drr, bend left=10, "pr_2"] && \\
U   & U\times \Delta^m \ar[d, "p_2"'] \ar[l, "p_1"]  & (U\times \Delta^m) \times (V \times \Delta^n) \ar[u, "\pi"] \ar[l, "\pi_1"'] \ar[d, "\rho"] \ar[r, "\pi_2"] & V \times \Delta^n \ar[d, "q_2"'] \ar[r, "q_1"']  &   V \\
  & \Delta^n  & \Delta^m \times \Delta^n\ar[l, "\rho_1"]   \ar[r, "\rho_2"'] & \Delta^n & 
 \end{tikzcd}
 \]
For simplicity, let us assume here that the currents \( R\)  and \( S \) satisfy the conditions in \cite[Prop. 2.2]{MR4498559}. For more general currents, the proof is similar and follows from the techniques and results in \cite[\S 2]{MR4498559}.  

Consider test forms \( \varphi \in \scrA^*_c(U) \) and \( \varrho\in \scrA^*_c(V) \). Then
\[
 \begin{alignedat}{4}
&    (A\times B)^\vee_{R\times S}&&  \{ pr_1^*\varphi \wedge pr_2^*\varrho \}  := 
 (A\times B)^\#(R \times S) \left\{ \pi^* (pr_1^*\varphi \wedge pr_2^*\varrho) \right\} &&  \\
 &  &&  = 
\left[  (\llb A \rrb \times \llb B \rrb) \cap ( \llb  U \rrb \times R) \times (\llb V \rrb) \times S) \right] \{\pi_1^*p_1^* \varphi \wedge \pi_2^*q_1^* \varrho\} && \quad\quad \text{ (by\ \  \cite[Prop. 2.2.a]{MR4498559})  } \\
& && = 
  \left[ \llb A \rrb \cap  ( \llb  U \rrb \times R) \right] \times \left[ \llb B \rrb   \cap  (\llb V \rrb \times S) \right]  \{\pi_1^*p_1^* \varphi \wedge \pi_2^*q_1^* \varrho\} && \quad\quad \text{ (by\ \  \cite[Thm. 5.8(10)]{MR315561})  }   \\
& &&   =   \left[ \llb A \rrb \cap  ( \llb  U \rrb \times R) \right] (\pi_1^*p_1^* \varphi) \cdot  \left[ \llb B \rrb   \cap  (\llb V \rrb \times S) \right] (\pi_2^*q_1^* \varrho ) && \\
  & && =
   \pi_{1\#}\left[ \llb A \rrb \cap  ( \llb  U \rrb \times R) \right] (p_1^* \varphi) \cdot  \pi_{2\#} \left[ \llb B \rrb   \cap  (\llb V \rrb \times S) \right] ( q_1^* \varrho ) && \\
& &&    =
   A^\vee_R  (p_1^* \varphi)  \cdot B^\vee_S( q_1^* \varrho )  \ =: \    (A^\vee_R  \times  B^\vee_S) (p_1^* \varphi)\wedge q_1^* \varrho  ), &&
   \end{alignedat}
   \]
where the last four identities follow from definitions and \cite[\S 4.1.8]{MR0257325}.
This proves the first assertion.

Now, let us use the commutative diagram diagram
\[
\begin{tikzcd}
U\times V & U\times V \times \Delta^{m+n}  \ar[l, "\hat \pi"']\ar[r,"\hat \rho"] \ar[d, "1 \times \lambda_\tau"', bend right=15] & \Delta^{m+n} \ar[d, "\lambda_\tau"]\\
& (U\times \Delta^m) \times (V\times \Delta^n) \ar[r, "\pi_2"'] \ar[u, "1\times \lambda^{-1}_\tau"', bend right=15] \ar[ul, "\pi"] & \Delta^m \times \Delta^n
\end{tikzcd}
\]
as follows. 
Given  a shuffle \( \tau \in \shuf{m}{n} \), let \( (1 \times \lambda_\tau) \colon U\times V \times \Delta^{m+n} \to (U\times \Delta^m)\times (V\times \Delta^n) \) be the isomorphism induced by \( \lambda_\tau \colon \Delta^{m+n} \xrightarrow{\ \cong \ } \Delta^m \times \Delta^n \), and observe that \( (1 \times \lambda^{-1}_\tau) = (1\times \lambda_\tau)^{-1}\) and hence \( (1 \times \lambda_\tau)^*(A\times B)=  (1 \times \lambda_\tau^{-1})_*(A\times B) \).
Therefore, the same arguments used  above give
\begin{align*}
& \{ (1 \times \lambda_\tau)^*(A\times B)\}^\vee_T  = 
\hat \pi_\# \{ (1 \times \lambda_\tau)^*(A\times B)\}^\#_T   = 
\hat \pi_\# \{ \llb (1 \times \lambda_\tau)^*(A\times B)\rrb \cap  ( \llb U \times V \rrb \times T)  \}\\
& =
\hat \pi_\# \{ \llb (1 \times \lambda_\tau^{-1})_*(A\times B)\rrb \cap  ( \llb U \times V \rrb \times T) \}   = 
\hat \pi_\# \left\{ (1 \times \lambda_\tau^{-1})_\# \left( \llb  A\times B\rrb \cap  (1\times \lambda_\tau)_\# ( \llb U \times V \rrb \times T)\right) \right\}\\
& = 
\{\hat \pi\circ  (1 \times \lambda_\tau^{-1})\}_\# \left\{\left( \llb  A\times B\rrb \cap  (1\times \lambda_\tau)_\# ( \llb U \times V \rrb \times T)\right) \right\} \\
& = 
\pi_\# \left\{\left( \llb  A\times B\rrb \cap ( \llb U \times V \rrb \times \lambda_{\tau \#}T)\right) \right\}  =
\pi_\#  \left(  A\times B\right)^\#_{\lambda_{\tau \#}T}   \ = \ (A\times B)^\vee_{\lambda_{\tau\#}T}.
\end{align*}

It follows that
\begin{align*}
(A\bar \times B)^\vee_T & = \sum_{\tau \in \shuf{m}{n}}\ \epsilon(\tau) \left\{ (1\times \lambda_\tau )^*(A\times B) \right\}^\vee_T
 = \sum_{\tau \in \shuf{m}{n}}\ \epsilon(\tau) (A\times B)^\vee_{\lambda_{\tau\#}T} \ = \ 
(A\times B)^\vee_{  \sum_{\tau }\epsilon(\tau) \lambda_{\tau\#}T}\\
& =
(A\times B)^\vee_{\psi_{m,n\#}(T)}.
\end{align*}
 \end{proof}

\bibliographystyle{amsplain}
%\bibliography{References.bib}

\providecommand{\bysame}{\leavevmode\hbox to3em{\hrulefill}\thinspace}
\providecommand{\MR}{\relax\ifhmode\unskip\space\fi MR }
\providecommand{\MRhref}[2]{%
  \href{http://www.ams.org/mathscinet-getitem?mr=#1}{#2}
}
\providecommand{\href}[2]{#2}

\end{document}